\documentclass[11pt]{article}

\usepackage{hyperref}

\usepackage{array} 
\usepackage{amsfonts}
\usepackage{amsmath}
\usepackage{amssymb}
\usepackage{graphicx}
\usepackage{float} 
\usepackage{color}
\usepackage{esint}    

\usepackage{enumitem} 

\newtheorem{thm}{Theorem}[section]
\newtheorem{cor}[thm]{Corollary}
\newtheorem{lem}[thm]{Lemma}
\newtheorem{prop}[thm]{Proposition}
\newtheorem{defn}[thm]{Definition}

\numberwithin{equation}{section}

\newcommand{\dx}{\,{\rm d}x}
\newcommand{\dy}{\,{\rm d}y}

\newcommand{\dt}{\,{\rm d}t}
\newcommand{\rd}{{\rm d}}

\def\LL{\mathrm{L}} 

\newcommand{\A}{\mathcal{L}}

\newcommand{\AI}{\mathcal{L}^{-1}}

\newcommand{\n}{F}

\newcommand{\p}{{\delta^\gamma}} 

\newcommand{\ka}{\overline{\kappa}}
\newcommand{\ua}{\overline{u}}
\newcommand{\kb}{\underline{\kappa}}
\newcommand{\ub}{\underline{u}}
\newcommand{\ud}{{u_\delta}}

\renewcommand{\k}{\kappa}
\newcommand{\K}{{\mathbb G}}

\newcommand{\RR}{\mathbb{R}}

\def\ee{\mathrm{e}} 
\def\dist{\mathrm{dist}} 
\def\dom{\mathrm{dom}} 

\def\qed{\,\unskip\kern 6pt \penalty 500
\raise -2pt\hbox{\vrule \vbox to8pt{\hrule width 6pt
\vfill\hrule}\vrule}\par}

\def\quotient#1#2{\raise1ex\hbox{$#1$}\Big/\lower1ex\hbox{$#2$}}

\definecolor{darkblue}{rgb}{0.05, .05, .65}
\definecolor{darkgreen}{rgb}{0.1, .65, .1}
\definecolor{darkred}{rgb}{0.8,0,0}

\begin{document}
\title{\textbf{Sharp global estimates 
for local and nonlocal \\ porous medium-type equations \\ in bounded domains}}

\author{\Large Matteo Bonforte, 
Alessio Figalli, 
~and~ Juan Luis V\'azquez 
}
\date{} 

\maketitle


\begin{abstract}
This paper provides a quantitative study of nonnegative solutions to  nonlinear diffusion equations
of porous medium-type of the form \ $\partial_t u + \A u^m=0$, $m>1$, where the operator $\A$ belongs to a general class of linear operators, and the equation is posed in a bounded domain $\Omega\subset \RR^N$.
As possible operators we  include the three most common definitions of the fractional Laplacian in a bounded domain with zero Dirichlet  conditions, and also a number of other nonlocal versions. In particular,  $\A$
can be a power of a uniformly elliptic operator with $C^1$ coefficients. Since
the nonlinearity is given by $u^m$
with $m>1$, the equation is degenerate parabolic.

The basic well-posedness theory for this class of equations has been recently developed in \cite{BV-PPR1,BV-PPR2-1}.
Here we address the regularity theory: decay and positivity, boundary behavior, Harnack inequalities, interior and boundary regularity, and asymptotic behavior. All this is done in a quantitative way, based on sharp a priori estimates. Although our focus is on the fractional models, our results cover also the local case when $\A$ is a uniformly elliptic operator, and provide new estimates even in this setting.

A surprising aspect discovered in this paper  is the possible presence of non-matching powers for the long-time boundary behavior. More precisely, when $\A=(-\Delta)^s$ is a spectral power of the {Dirichlet} Laplacian inside a smooth domain,
we can prove that: \\
- when $2s> 1-1/m$,  for large times all solutions behave as ${\rm dist}^{1/m}$ near the boundary; \\
- when $2s\le 1-1/m$,  different solutions may exhibit different boundary behavior.\\
This unexpected phenomenon is a completely new feature of the nonlocal nonlinear structure of this model, and it is not present in the semilinear elliptic equation $\A u^m=u$.

\end{abstract}

\vskip 1 cm


\noindent {\sc Addresses:}

\noindent Matteo Bonforte. Departamento de Matem\'{a}ticas, Universidad
Aut\'{o}noma de Madrid,\\ Campus de Cantoblanco, 28049 Madrid, Spain.
e-mail address:~\texttt{matteo.bonforte@uam.es}

\noindent Alessio Figalli. ETH Z\"urich, Department of Mathematics,
R\"amistrasse 101,\\ 8092 Z\"urich, Switzerland.
E-mail:\texttt{~alessio.figalli@math.ethz.ch}

\noindent Juan Luis V\'azquez. Departamento de Matem\'{a}ticas, Universidad
Aut\'{o}noma de Madrid,\\ Campus de Cantoblanco, 28049 Madrid, Spain.  e-mail address:~\texttt{juanluis.vazquez@uam.es}
 
\bigskip


\medskip


\newpage
\footnotesize
\tableofcontents
\normalsize

\newpage
\section{Introduction}

In this paper we address the question of obtaining a priori estimates, positivity, boundary behavior, Harnack inequalities, and regularity for a suitable class of  weak solutions of  nonlinear nonlocal diffusion equations of the form:\vspace{-2mm}
\begin{equation}\label{FPME.equation}
\partial_t u+\A\,\n (u)=0  \qquad\mbox{posed in }Q=(0,\infty)\times \Omega\,,\vspace{-2mm}
\end{equation}
where $\Omega\subset \RR^N$ is a bounded domain with $C^{1,1}$ boundary, $N\ge 2$\footnote{Our results work also in dimension $N=1$ if the fractional exponent (that we shall introduce later) belongs to the range $0<s<1/2$. The interval $1/2\le s<1$ requires some minor modifications that we prefer to avoid in this paper.}, and $\A$ is a linear operator representing diffusion of local or nonlocal type, the prototype example being the fractional Laplacian (the class of admissible operators will be precisely described below).
Although our arguments hold for a rather general class of nonlinearities  $F:\RR \to \RR$, for the sake of simplicity we shall focus on the model case $F(u)=u^m$ with $m>1$.

The use of nonlocal operators in diffusion equations reflects the need to model the presence of long-distance effects not included in evolution driven by the Laplace operator, and this is well documented in the literature. The physical motivation and relevance of the nonlinear diffusion models with nonlocal operators has been mentioned in many references, see for instance \cite{AC,BV2012,BV-PPR1,DPQRV1,DPQRV2,Vaz2014}.
Because $u$ usually represents a density, all data and solutions are supposed to be nonnegative. Since the problem is posed in a bounded domain we need boundary or external conditions that we assume of Dirichlet type.

This kind of problems has been extensively studied when $\A=-\Delta$ and  $\n(u)=u^m$, $m>1$,
in which case the equation becomes the classical Porous Medium Equation \cite{VazBook,DK0, DaskaBook, JLVmonats}. 
Here, we are interested in treating nonlocal diffusion operators, in particular fractional Laplacian operators.  Note that, since we are working on a bounded domain,  the concept of fractional Laplacian operator admits  several non-equivalent versions, the best known being the Restricted   Fractional Laplacian (RFL), the Spectral Fractional Laplacian (SFL), and the Censored  Fractional Laplacian (CFL);
see Section \ref{ssec.examples} for more details. We use these names because they already appeared in some previous works \cite{BSV2013, BV-PPR2-1}, but we point out that RFL is usually known as the Standard Fractional Laplacian, or plainly Fractional Laplacian, and CFL is often called Regional Fractional Laplacian.   

The case of the SFL operator with $F(u)=u^m$, $m>1$,  has been already studied by the first and the third author in  \cite{BV-PPR1,BV-PPR2-1}. In particular, in \cite{BV-PPR2-1} the authors presented a rather abstract setting where they were able to treat not only the usual fractional Laplacians but also a large number of variants that will be listed below for the reader's convenience. Besides, rather general increasing nonlinearities $F$ were allowed. The basic questions of existence and uniqueness of suitable solutions for this problem were solved in \cite{BV-PPR2-1} in the class of `weak dual solutions', an extesion of the concept of solution introduced in \cite{BV-PPR1} that has proved to be quite flexible and efficient. A number of a priori estimates (absolute bounds and smoothing effects) were also derived in that generality.

Since these basic facts are settled, here we focus our attention on  the finer aspects of the theory,  mainly  sharp boundary estimates and decay estimates.  Such upper and lower bounds will be formulated in terms of the first eigenfunction $\Phi_1$ of $\A$, that under our assumptions will satisfy $\Phi_1\asymp  \dist(\cdot, \partial\Omega)^\gamma$ for a certain characteristic power $\gamma\in (0,1]$ that depends on the particular operator we consider. Typical values are $\gamma=s$ (SFL), $\gamma=1$ (RFL), and $\gamma=s-1/2$ for $s>1/2$ (CFL), cf. Subsections \ref{ssec.examples} and \ref{sec.examples}. As a consequence, we get various kinds of local and global Harnack type inequalities.

It is worth mentioning that  some of the boundary estimates that we obtain  for the parabolic case are essentially elliptic in nature. The study of this issue for stationary problems is done in a companion paper \cite{BFV-Elliptic}. This has the advantage that many arguments are clearer, since the parabolic problem is more complicated than the elliptic one.  Clarifying  such difference is one of the main contributions of our present work.

Thanks to these results, in the last part of the paper we are  able to prove  both interior and boundary regularity, and to find the large-time asymptotic behavior of solutions.

Let us indicate here some notation of general use. The symbol $\infty$ will always denote $+\infty$. Given $a,b$, we use the notation $a\asymp b$ whenever there exist universal constants $c_0,c_1>0$ such that $c_0\,b\le a\le c_1 b$\,. We also use the symbols $a\vee b=\max\{a,b\}$ and $a\wedge b=\min\{a,b\}$.
We will always consider bounded domains $\Omega$ with boundary of class $C^{2}$. In the paper we use the short form `solution' to mean `weak dual solution', unless differently stated.
%
%
\subsection{Presentation of the  results on sharp boundary behaviour}\label{ssec.results.boundary}
\noindent $\bullet$ A basic principle in the paper is that the sharp boundary estimates depend not only on $\A$ but also on the behavior  of the nonlinearity $\n(u)$ near $u=0$, i.e., in our case, on the exponent $m>1$.  The elliptic analysis performed in the companion paper \cite{BFV-Elliptic} combined with some standard arguments will allow us to prove that, in {\em all} cases, $u(t)$ approaches the separate-variable solution ${\mathcal U}(x,t)=t^{-\frac1{m-1}}S(x)$ in the sense that
\begin{equation}\label{asymp.intro}
\left\|t^{\frac{1}{m-1}}u(t,\cdot)- S\right\|_{\LL^\infty(\Omega)}\xrightarrow{t\to\infty}0,
\end{equation}
where $S$ is the solution of the elliptic problem
(see Theorems \ref{Thm.Elliptic.Harnack.m} and \ref{Thm.Asympt.0}).
The behaviour of the profile $ S(x)$ is shown to be,  when $2sm\ne \gamma(m-1)$,
\begin{equation}\label{as.sep.var}
S(x)\asymp \Phi_1(x)^{\sigma/m},  \qquad \sigma:=\min\left\{1,\frac{2sm}{\gamma(m-1)}\right\}.
\end{equation}
Thus, the behavior strongly depends on the new parameter $\sigma$, more precisely, on whether this parameter is equal to $1$ or less than $1$. As we shall see later, $\sigma$ encodes the interplay between the ``elliptic scaling power'' $2s/(m-1)$, the ``eigenfunction power'' $\gamma$, and the ``nonlinearity power'' $m$. When $2sm =\gamma(m-1)$ we have $\sigma=1$, but a logarithmic correction appears:
\begin{equation}\label{as.sep.var.1}
S(x)\asymp \Phi_1(x)^{1/m}\left(1+|\log\Phi_1(x) |\right)^{1/(m-1)}\,.
\end{equation}
\noindent $\bullet$  This fact and the results in \cite{BFR} prompted us to look for estimates of the form
\begin{equation}\label{intro.1b}
c_0(t)\frac{\Phi_1^{\sigma/m}(x_0)} {t^{\frac1{m-1}}}
\le u(t,x_0) \le c_1\frac{ \Phi_1^{\sigma/m}(x_0)} {t^{\frac1{m-1}}}\qquad \text{for all $t>0$, $x_0\in\Omega$,}
\end{equation}
where $c_0(t)$ and $c_1$ are positive and independent of $u$, eventually with a logarithmic term appearing when $2sm =\gamma(m-1)$, as in \eqref{as.sep.var.1}.  We will prove in this paper that the upper bound holds for the three mentioned Fractional Laplacian choices,  and indeed for the whole class of integro-differential operators we will introduce below, cf.  Theorem \ref{thm.Upper.PME.II}. Also, separate-variable solutions saturate the upper bound.

The issue of the validity of a lower bound as in \eqref{intro.1b} is instead much more elusive.
A first indication for this is the introduction of a function $c_0(t)$ depending on $t$, instead of a constant. This wants to reflect the fact that the solution may take some time to reach the boundary behaviour that is expected to hold uniformly for large times. Indeed, recall that in the classical PME \cite{Ar-Pe,JLVmonats, VazBook}, for data supported away from the boundary, some `waiting time' is needed for the support to reach the boundary.

\medskip

\noindent $\bullet$ As proved in \cite{BFR}, the stated lower bound holds for the RFL with $c_0(t)\sim (1\wedge t)^{m/(m-1)}.$
In particular, in this nonlocal setting, infinite speed of propagation holds.
Here, we show that this holds also for the CFL  and a number of other operators, cf.  Theorem \ref{Thm.lower.B}.
Note that for the RFL and CFL we have  $2sm>\gamma(m-1)$, in particular $\sigma=1$  which simplifies formula \eqref{intro.1b}.\\ A combination of an upper and a lower bound with matching behaviour (with respect to $x$ and $t$) will be called a {\sl Global Harnack Principle}, and holds for all  $t>0$ for these operators, cf.  Theorems \ref{thm.GHP.PME.I} and \ref{thm.GHP.PME.II}.

\medskip

\noindent $\bullet$ When  $\A$ is the SFL, we shall see that the lower bound may fail. Of course,  solutions by separation of variables satisfy the matching estimates in \eqref{intro.1b} (eventually with a extra logarithmic term in the limit case, as in \eqref{as.sep.var.1}),  but it came as a complete surprise to us that for the SFL the situation is not the same for ``small'' initial data. More precisely:

\noindent(i) We can  prove that the following  bounds always hold for all times:
\begin{equation}\label{intro.1}
c_0\left(1\wedge \frac{t}{t_*}\right)^{\frac{m}{m-1}}\frac{\Phi_1(x_0)} {t^{\frac1{m-1}}}
\le u(t,x_0) \le c_1\frac{ \Phi_1^{\sigma/m}(x_0)} {t^{\frac1{m-1}}}\,,
\end{equation}
(when $2sm=\gamma(m-1)$, a logarithmic correction $\left(1+|\log\Phi_1(x) |\right)^{1/(m-1)}$ appears in the right hand side),
cf.  Theorem \ref{thm.Lower.PME}.
 These are non-matching estimates.

\noindent(ii) For $2sm>\gamma(m-1)$, the sharp estimate \eqref{intro.1b} holds for any nonnegative nontrivial solution  for large times  $t \geq t_*$, cf.  Theorem \ref{thm.Lower.PME.large.t}.

\noindent(iii) \textbf{Anomalous boundary behaviour.} Consider now the SFL with $\sigma<1$ (resp.  $2sm=\gamma(m-1)$).\footnote{Since for the SFL $\gamma=1$, we have $\sigma<1$ if and only if
$$
0<s<s_*:=\frac{m-1}{2m}<\frac{1}{2}.
$$
Note that $s_*\to 0$ as we tend to the linear case $m=1$, so this exceptional regime dooes not appear for linear diffusions, both fractional and standard.}
In this case we can find initial data for which the upper bound in \eqref{intro.1} is not sharp. Depending on the initial data, there are several possible rates for the long-time behavior near the boundary. More precisely:
\begin{enumerate}
\item[(a)] When $u_0\leq A\,\Phi_1$, then $u(t)\le F(t)\Phi_1^{1/m} \ll \Phi_1^{\sigma/m}$ (resp. $\Phi_1^{1/m} \ll \Phi_1^{1/m}\left(1+|\log\Phi_1 |\right)^{1/(m-1)}$) for all times, see Theorem \ref{prop.counterex}. In particular
\begin{equation}\label{limit.intro}
\lim_{x\to \partial\Omega}\frac{u(t,x)}{\Phi_1(x)^{\sigma/m}}= 0  \quad \Bigl(\text{resp.} \lim_{x\to \partial\Omega}\frac{u(t,x)}{\Phi_1(x)^{1/m}\left(1+|\log\Phi_1(x) |\right)^{1/(m-1)}}= 0\Bigr)
\end{equation}
for any $t>0.$
\item[(b)] When  $u_0\leq A\,\Phi_1^{1-2s/\gamma}$ then $u(t)\le F(t)\Phi_1^{1-2s/\gamma}$ for small times, see Theorem \ref{thm.Upper.PME.III}. Notice that when $\sigma<1$ we have always $1-\frac{2s}{\gamma}>\sigma/m$.
This sets a limitation to the improvement of the lower bound, which is confirmed by another result: In Theorem \ref{prop.counterex2} we show that lower bounds of the form $u(T,x)\geq \kb\Phi_1^\alpha(x)$ for  data $u_0(x)\le A\Phi_1(x)$ are possible only for $\alpha\ge 1-2s/\gamma$.
\item[(c)] On the other hand, for ``large'' initial data, Theorem \ref{thm.GHP.PME.II} shows  that the desired matching estimates from above and below hold.
\end{enumerate}
 After discovering this strange boundary behavior, we looked for  numerical confirmation. In Section \ref{sec.numer} we will explain the numerical results obtained in \cite{numerics}.  Note that, if one looks for universal bounds independent of the initial condition, Figures 2-3 below seem to suggest that the bounds provided by \eqref{intro.1} are optimal for all times and all operators.

\medskip

\noindent $\bullet$  The current interest in more general types of nonlocal operators led us to a more general analysis where the just explained alternative has been extended to a wide class of integro-differential operators, subject only to a list of properties that we call (A1), (A2), (L1), (L2), (K2), (K4); a number of examples are explained in Section \ref{sec.hyp.L}. These general classes appear also in the study of the elliptic problem \cite{BFV-Elliptic}.


\subsection{Asymptotic behaviour and regularity}

Our quantitative lower and upper estimates admit a formulation as local or global Harnack inequalities. They are used at the end of the paper to settle two important issues.

\noindent\textbf{Sharp asymptotic behavior. }Exploiting the techniques in \cite{BSV2013},
we can prove a sharp asymptotic behavior for our nonnegative and nontrivial solutions when the upper and lower bound have matching powers. Such sharp results hold true for a quite general class of local and nonlocal operators. A detailed account is given in Section \ref{sec.asymptotic}.

\smallskip

\noindent\textbf{Regularity. }By a variant of the techniques used in \cite{BFR},
we can show interior H\"older regularity.
In addition, if the kernel of the operator satisfies some suitable continuity assumptions, we show that solutions are classical in the interior and are H\"older continuous up to the boundary if the upper and lower bound have matching powers.
We refer to Section \ref{sect.regularity} for  details.

\section{General class of operators and their kernels}\label{sec.hyp.L}

The interest of the theory developed here lies both in the sharpness of the results and in the wide range of applicability. We have just mentioned the most relevant examples appearing in the literature, and more are listed at the end of this section. Actually, our theory  applies to a general class of operators with definite assumptions, and this is what we want to explain now.

Let us present the properties that have to be assumed on the class of admissible operators.
Some of them already appeared in \cite{BV-PPR2-1}.
However, to further develop our theory, more hypotheses need to be introduced.
In particular, while \cite{BV-PPR2-1} only uses the properties of the Green function, here we shall make some assumptions also on the kernel of $\A$ (whenever it exists). Note that assumptions on the kernel $K$ of $\A$
 are needed for the positivity results, because we need to distinguish between the local and nonlocal cases. The study of the kernel $K$ is performed in Subsection \ref{ss2.2}.

For convenience of reference, the list of used assumptions is (A1), (A2), (K2), (K4), (L1), (L2).
The first three are assumed in all operators $\A$ that we use.

\medskip

\noindent $\bullet$ {\bf Basic assumptions on $\A$.}
The linear operator $\A: \dom(\A)\subseteq\LL^1(\Omega)\to\LL^1(\Omega)$ is assumed to be densely defined  and sub-Markovian, more precisely, it satisfies (A1) and (A2) below:
\begin{enumerate}
\item[(A1)] $\A$ is $m$-accretive on $\LL^1(\Omega)$;
\item[(A2)] If $0\le f\le 1$ then $0\le \ee^{-t\A}f\le 1$.
\end{enumerate}
Under these assumption, in \cite{BV-PPR2-1}, the first and the third author proved existence, uniqueness, weighted estimates, and smoothing effects.

\medskip

\noindent $\bullet$ {\bf Assumptions on the kernel.}
Whenever $\A$ is defined in terms of a kernel $K(x,y)$ via the formula
\begin{equation*} 
\A f(x)=P.V.\int_{\RR^N} \big(f(x)-f(y)\big)\,K(x,y)\dy\,,
\end{equation*}
assumption (L1) states that there exists $\underline{\kappa}_\Omega>0$ such that
\[\tag{L1}
\inf_{x,y\in \Omega}K(x,y)\ge \underline{\kappa}_\Omega>0\,.
\]
We note that condition holds both for  the RFL and the CFL, see Section \ref{ssec.examples}.

\noindent- Whenever $\A$ is defined in terms of a kernel $K(x,y)$ and a zero order term via the formula
\[
\A f(x)=P.V.\int_{\RR^N} \big(f(x)-f(y)\big)\,K(x,y)\dy + B(x)f(x),
\]
assumptions (L2) states that
\[\tag{L2}
K(x,y)\ge c_0\p(x)\p(y),\quad c_0>0, \qquad\mbox{and}\quad B(x)\ge 0,
\]
where, from now on, we adopt the notation $\delta(x):=\dist(x, \partial\Omega)$.
This condition is satisfied by the SFL in a stronger form, see Section \ref{ss2.2} and Lemma \ref{Lem.Spec.Ker}.

\medskip

\noindent $\bullet$ {\bf Assumptions on $\AI$.} In order to prove our quantitative estimates, we need to be more specific about the operator $\A$. Besides satisfying (A1) and (A2), we will assume that it has a left-inverse $\AI: \LL^1(\Omega)\to \LL^1(\Omega)$  that can be represented by a kernel $\K$ (the letter ``G'' standing for Green function) as
\[
\AI[f](x)=\int_\Omega \K(x,y)f(y)\dy\,,
\]
where $\K$ satisfies the following assumption, for some $s\in (0,1]$:
There exist constants $\gamma\in (0,1]$ and $c_{0,\Omega},c_{1,\Omega}>0$ such that, for a.e. $x,y\in \Omega$,
\[\tag{K2}
c_{0,\Omega}\,\p(x)\,\p(y) \le \K(x,y)\le \frac{c_{1,\Omega}}{|x-y|^{N-2s}}
\left(\frac{\p(x)}{|x-y|^\gamma}\wedge 1\right)
\left(\frac{\p(y)}{|x-y|^\gamma}\wedge 1\right).
\]
(Here and below we use the labels (K2) and (K4) to be consistent with the notation in \cite{BV-PPR2-1}.)
Hypothesis (K2) introduces an exponent $\gamma$ which is a  characteristic of the operator and will play a big role in the results.
Notice that defining an inverse operator $\AI$ implies that we are taking into account the Dirichlet boundary conditions. See more details in Section 2 of  \cite{BV-PPR2-1}.

\medskip

\noindent - The lower bound in (K2) is weaker than the known bounds on the Green function for many examples under consideration; indeed, the following stronger estimate holds in many cases:
\[\tag{K4}
\K(x,y)\asymp \frac{1}{|x-y|^{N-2s}}
\left(\frac{\p(x)}{|x-y|^\gamma}\wedge 1\right)
\left(\frac{\p(y)}{|x-y|^\gamma}\wedge 1\right)\,.
\]
\noindent\textbf{Remarks. }(i) The labels (A1), (A2), (K1), (K2), (K4) are consistent with the notation in \cite{BV-PPR2-1}. The label (K3) was used to mean hypothesis (K2) written in terms of $\Phi_1$ instead of $\p$.\\
(ii) In the classical local case $\A=-\Delta$, the Green function $\K$ satisfies (K4) only when $N\geq 3$, as the formulas slightly change when $N=1,2$. In the fractional case $s \in (0,1)$
the same problem arises when $N=1$ and $s \in [1/2,1)$. Hence, treating also these cases would require a slightly different analysis based on different but related assumptions on $\K$. Since our approach is very general, we expect it to work also in these remaining cases without any major difficulties. However, to simplify the presentation, from now on we assume that
$$
\text{either $N\geq 2$ and $s\in(0,1),\qquad$ or $N=1$ and $s \in (0,1/2)$.}
$$

\noindent\textbf{The role of the first eigenfunction of $\A$. }We have shown in \cite{BFV-Elliptic} that, under assumption (K1), the operator $\A$ is compact, it has a discrete spectrum, and a first nonnegative bounded eigenfunction $\Phi_1$; assuming also (K2), we have that
\begin{equation}\label{Phi1.est}
\Phi_1(x)\asymp \p(x)=\dist(x,\partial\Omega)^\gamma\qquad\mbox{for all }x\in \overline{\Omega}.
\end{equation}
Hence, $\Phi_1$ encodes the parameter $\gamma$ that takes care of describing the boundary behavior. We recall that we are assuming that the boundary of $\Omega$ is smooth enough, for instance $C^{1,1}$.

\noindent\textbf{Remark. }We note that our assumptions allow us to cover all the examples of operators described in Sections \ref{ssec.examples} and \ref{sec.examples}.

\subsection{Main examples of operators and properties}
\label{ssec.examples}

When working in the whole $\mathbb R^N$, the fractional Laplacian admits different definitions that can be shown to be all equivalent. On the other hand, when we deal with bounded domains, there are at least three different operators in the literature, that we call the Restricted (RFL), the Spectral (SFL)  and the Censored Fractional Laplacian (CFL). We will show below that these different operators exhibit quite different behaviour, so the distinction between them has to be taken into account.
Let us present the statement and results for the three model cases, and we refer to Section \ref{sec.examples} for further examples.  Here, we collect the sharp results about the boundary behavior, namely the Global Harnack inequalities from Theorems \ref{thm.GHP.PME.I}, \ref{thm.GHP.PME.II}, and \ref{thm.GHP.PME.III}.

\medskip

\noindent\textit{The parameters $\gamma$ and $\sigma$.} The strong difference between the various operators $\A$ is reflected in the different boundary behavior of their nonnegative solutions. We will often use the exponent $\gamma$, that represents the boundary behavior of the first eigenfunction $\Phi_1 \asymp \dist(\cdot,\partial\Omega)^\gamma$, see~\cite{BFV-Elliptic}.
Both in the parabolic theory of this paper and the elliptic theory of paper \cite{BFV-Elliptic}
the parameter $\sigma=\min\left\{1, \frac{2sm}{\gamma(m-1)} \right\}$ introduced in    \eqref{as.sep.var} plays a big role.
\subsubsection{The RFL} We define the fractional Laplacian operator acting on a bounded domain by using the integral representation on the whole space in terms of a hypersingular kernel, namely
\begin{equation}\label{sLapl.Rd.Kernel}
(-\Delta_{\RR^N})^{s}  g(x)= c_{N,s}\mbox{
P.V.}\int_{\mathbb{R}^N} \frac{g(x)-g(z)}{|x-z|^{N+2s}}\,dz,
\end{equation}
where $c_{N,s}>0$ is a normalization constant, and we ``restrict'' the operator to functions that are zero outside $\Omega$. We denote such operator by $\A=(-\Delta_{|\Omega})^s$\,, and call it the \textit{restricted fractional Laplacian}\footnote{In the literature this is often called the fractional Laplacian on domains, but this simpler name may be confusing when the spectral fractional Laplacian is also considered, cf. \cite{BV-PPR1}. As discussed in this paper, there are other natural versions.} (RFL). The initial and boundary conditions associated to the fractional diffusion equation \eqref{FPME.equation} read
$u(t,x)=0$ in $(0,\infty)\times\RR^N\setminus \Omega$ and $u(0,\cdot)=u_0$. As explained in \cite{BSV2013}, such boundary conditions can also be understood via the Caffarelli-Silvestre extension, see \cite{Caffarelli-Silvestre}. The sharp expression of the boundary behavior for RFL has been investigated in \cite{RosSer}. We refer to \cite{BSV2013} for a careful construction of the RFL in the framework of fractional Sobolev spaces, and \cite{BlGe} for a probabilistic interpretation.

This operator satisfies the assumptions (A1), (A2), (L1),  and also (K2) and (K4) with $\gamma=s<1$. Let us present our  results in this case. Note that we have $\sigma=1$ for all $0<s < 1$, and Theorem \ref{thm.GHP.PME.I} shows the sharp boundary behavior for all times, namely for all $t>0$ and a.e. $x\in \Omega$ we have
\begin{equation}\label{thm.GHP.PME.I.Ineq.00}
\kb\, \left(1\wedge \frac{t}{t_*}\right)^{\frac{m}{m-1}}\frac{\dist(x,\partial\Omega)^{s/m}}{t^{\frac{1}{m-1}}}
\le \, u(t,x) \le \ka\, \frac{\dist(x,\partial\Omega)^{s/m}}{t^{\frac{1}{m-1}}}\,.
\end{equation}
The critical time $t_*$ is given by a weighted $\LL^1$ norm, namely $t_*:= \k_* \|u_0\|_{\LL^1_{\Phi_1}(\Omega)}^{-(m-1)}$,
where $\k_*>0$ is a universal constant. Moreover,  solutions are classical in the interior and we prove sharp H\"older continuity up to the boundary.
These regularity results have been first obtained in \cite{BFR}; we give here different proofs valid in the more general setting of this paper.  See Section \ref{sect.regularity} for further details.
\subsubsection{The SFL } Starting  from the classical Dirichlet Laplacian $\Delta_{\Omega}$ on the domain $\Omega$\,, the so-called {\em spectral definition} of the fractional power of $\Delta_{\Omega}$ may be defined via a formula in terms of the semigroup associated to the Laplacian, namely
\begin{equation}\label{sLapl.Omega.Spectral}
\displaystyle(-\Delta_{\Omega})^{s}
g(x)= \frac1{\Gamma(-s)}\int_0^\infty
\left(e^{t\Delta_{\Omega}}g(x)-g(x)\right)\frac{dt}{t^{1+s}}=\sum_{j=1}^{\infty}\lambda_j^s\, \hat{g}_j\, \varphi_j(x)\,,
\end{equation}
where  $(\lambda_j,\varphi_j)$, $j=1,2,\ldots$, is the normalized spectral sequence of the standard Dirichlet Laplacian on $\Omega$\,, $
\hat{g}_j=\int_\Omega g(x)\varphi_j(x)\dx$, and $\|\varphi_j\|_{\LL^2(\Omega)}=1$\,.
We denote this operator by $\A=(-\Delta_{\Omega})^s$\,, and call it the \textit{spectral fractional Laplacian} (SFL) as in \cite{Cabre-Tan}.  The initial and boundary conditions associated to the fractional diffusion equation \eqref{FPME.equation} read
$u(t,x)=0$ on $(0,\infty)\times\partial\Omega$ and $u(0,\cdot)=u_0$. Such boundary conditions can also be understood via the Caffarelli-Silvestre extension, see \cite{BSV2013}. Following ideas of \cite{SV2003}, we use the fact that this operator admits a kernel representation,
\begin{equation}\label{SFL.Kernel}
(-\Delta_{\Omega})^{s}  g(x)= c_{N,s}\mbox{
P.V.}\int_{\Omega} \left[g(x)-g(z)\right]K(x,z)\,dz + B(x)g(x)\,,
\end{equation}
where $K$ is a singular and compactly supported kernel, which degenerates at the boundary, and $B\asymp \dist(\cdot,\partial\Omega)^{-2s}$ (see \cite{SV2003} or Lemma \ref{Lem.Spec.Ker} for further details).
This operator satisfies the assumptions (A1), (A2), (L2),  and also (K2) and (K4) with $\gamma=1$. Therefore, $\sigma$ can be less than $1$, depending on the values of $s$ and $m$.

As we shall see, in our parabolic setting, the degeneracy of the kernel is responsible for a peculiar change of the boundary behavior of the solutions (with respect to the previous case) for small and large times. Here, the lower bounds change both for short and large times,
and they strongly depend on $\sigma$ and on $u_0$: we called this phenomenon \textit{anomalous boundary behaviour }in Subsection \ref{ssec.results.boundary}. More precisely, Theorem \ref{thm.GHP.PME.III} shows that for all $t>0$ and all $x\in \Omega$ we have

\begin{equation}\label{thm.GHP.PME.III.Ineq.00}
\kb\,\left(1\wedge \frac{t}{t_*}\right)^{\frac{m}{m-1}}\frac{\dist(x,\partial\Omega) }{t^{\frac{1}{m-1}}}
\le \, u(t,x) \le \ka\, \frac{\dist(x,\partial\Omega)^{\sigma/m}}{t^{\frac1{m-1}}}
\end{equation}
(when $2sm=\gamma(m-1)$, a logarithmic correction $\left(1+|\log\Phi_1(x) |\right)^{1/(m-1)}$ appears in the right hand side).
Such lower behavior is somehow minimal, in the sense that it holds in all cases. The basic asymptotic result (cf. \eqref{asymp.intro} or Theorem \ref{Thm.Asympt.0.1}) suggests that the lower bound in \eqref{thm.GHP.PME.III.Ineq.00} could be improved by replacing $\dist(x,\partial\Omega)$ with $\dist(x,\partial\Omega)^{\sigma/m}$, at least for large times. This is shown to be true for $\sigma=1$ (cf. Theorem \ref{thm.Lower.PME.large.t}), but it is false for $\sigma<1$ (cf. Theorem \ref{prop.counterex}), since there are ``small'' solutions with non-matching boundary behaviour for all times, cf. \eqref{limit.intro}.

It is interesting that, in this case, one can appreciate the interplay between the ``elliptic scaling power'' $2s/(m-1)$ related to the invariance of the equation $\A S^m=S$ under the scaling $S(x)\mapsto \lambda^{-2s/(m-1)}S(\lambda x)$, the ``eigenfunction power'' $\gamma=1$,
and the ``nonlinearity power'' $m$, made clear through the parameter $\sigma/m$.
Also in this case, thanks to the strict positivity in the interior, we can show interior space-time regularity of solutions, as well as sharp boundary H\"older regularity for large times whenever upper and lower bounds match.
 
\subsubsection{The CFL} In the simplest case, the infinitesimal operator of the censored stochastic processes  has the form
\begin{equation}
\A g(x)=\mathrm{P.V.}\int_{\Omega}\frac{g(x)-g(y)}{|x-y|^{N+2s}}\dy\,,\qquad\mbox{with }\frac{1}{2}<s<1\,.
\end{equation}
This operator has been  introduced in \cite{bogdan-censor} (see also \cite{Song-coeff} and \cite{BV-PPR2-1} for further details and references).

In this case $\gamma=s-1/2<2s$, hence $\sigma=1$ for all $1/2<s < 1$,
and Theorem \ref{thm.GHP.PME.I} shows that for  all $t>0$ and $x\in \Omega$ we have
\[
\kb\, \left(1\wedge \frac{t}{t_*}\right)^{\frac{m}{m-1}}\frac{\dist(x,\partial\Omega)^{(s-1/2)/m}}{t^{\frac{1}{m-1}}}
\le \, u(t,x) \le \ka\, \frac{\dist(x,\partial\Omega)^{(s-1/2)/m}}{t^{\frac{1}{m-1}}}\,.
\]
Again, we have interior space-time regularity of solutions, as well as sharp boundary H\"older regularity for all times.

\medskip

\subsubsection{Other examples}
There a number of examples to which our theory applies, besides the RFL, CFL and SFL, since they satisfy the list of assumptions listed in the previous section.  Some are listed in the last Section \ref{sec.comm}, see more detail in \cite{BV-PPR2-1}.

%
\section{Reminders about weak dual solutions}\label{sec.results}

\medskip
We denote by $\LL^p_{\Phi_1}(\Omega)$ the weighted $\LL^p$ space $\LL^p(\Omega\,,\, \Phi_1\dx)$, endowed with the norm
\[
\|f\|_{\LL^p_{\Phi_1}(\Omega)}=\left(\int_{\Omega} |f(x)|^p\Phi_1(x)\dx\right)^{\frac{1}{p}}\,.
\]
 
\medskip

\noindent{\bf Weak dual solutions: existence and uniqueness.}
%
We recall the definition of weak dual solutions used in \cite{BV-PPR2-1}. This is expressed in terms of the inverse  operator $\AI$, and encodes the Dirichlet boundary condition.  This is needed to build a theory of bounded nonnegative unique solutions to Equation \eqref{FPME.equation} under the assumptions of the previous section. Note that in \cite{BV-PPR2-1} we have used the setup with the weight $\p=\dist(\cdot,\partial\Omega)^\gamma$, but the same arguments generalize immediately to the weight $\Phi_1$; indeed under assumption (K2), these two setups are equivalent.

\begin{defn}\label{Def.Very.Weak.Sol.Dual} A function $u$ is a {\sl weak dual} solution to the Dirichlet Problem for Equation \eqref{FPME.equation} in $(0,\infty)\times \Omega$ if:
\begin{itemize}[leftmargin=*]
\item $u\in C((0,\infty): \LL^1_{\Phi_1}(\Omega))$\,, $u^m \in \LL^1\left((0,\infty):\LL^1_{\Phi_1}(\Omega)\right)$;
\item  The identity
\begin{equation}
\displaystyle \int_0^\infty\int_{\Omega}\AI u \,\dfrac{\partial \psi}{\partial t}\,\dx\dt
-\int_0^\infty\int_{\Omega} u^m\,\psi\,\dx \dt=0
\end{equation}
holds for every test function $\psi$ such that  $\psi/\Phi_1\in C^1_c((0,\infty): \LL^\infty(\Omega))$\,.
\item A {weak dual} solution to the Cauchy-Dirichlet problem {(CDP)} is a  weak dual solution to the homogeneous Dirichlet Problem for equation \eqref{FPME.equation} such that $u\in C([0,\infty): \LL^1_{\Phi_1}(\Omega))$  and $u(0,x)=u_0\in \LL^1_{\Phi_1}(\Omega)$.
\end{itemize}\end{defn}
This kind of solution has been first introduced in \cite{BV-PPR1}, cf. also \cite{BV-PPR2-1}.
Roughly speaking, we are considering the weak solution to the ``dual equation'' $\partial_t U=- u^m$\,, where $U=\AI u$\,, posed on the bounded domain $\Omega$ with homogeneous Dirichlet conditions. Such weak solution is obtained by approximation from below as the limit of the unique mild solution provided by the semigroup theory (cf. \cite{BV-PPR2-1}),
and it was used in \cite{Vaz2012} with space domain $\RR^N$ in the study of Barenblatt solutions. We call those solutions \textit{minimal weak dual solutions, }and it has been proven in Theorems 4.4 and 4.5 of \cite{BV-PPR2-1} that such solutions exist and are unique for any nonnegative data $u_0\in\LL^1_{\Phi_1}(\Omega)$. The class of weak dual solutions includes the classes of weak, mild and strong solutions, and is included in the class of very weak solutions. In this class of solutions the standard comparison result holds.
 
\medskip

\noindent {\bf Explicit solution.} When trying to understand the  behavior of positive solutions with general nonnegative data, it is natural to look for solutions obtained by separation of variables. These are given by\vspace{-1mm}
\begin{equation}\label{friendly.giant}
\mathcal{U}_T(t,x):=(T+t)^{-\frac{1}{m-1}}S(x)\,,\qquad T\geq 0,\vspace{-1mm}
\end{equation}
where $S$ solves the elliptic problem\vspace{-1mm}
\begin{equation}\label{Elliptic.prob}
\left\{\begin{array}{lll}
\A S^m=  S &  ~ {\rm in}~ (0,+\infty)\times \Omega,\\
S=0 & ~\mbox{on the  boundary.}
\end{array}
\right.\vspace{-1mm}
\end{equation}
The properties of $S$ have been thoroughly studied in the companion paper \cite{BFV-Elliptic}, and  we summarize them here for the reader's convenience.\vspace{-2mm}
\begin{thm}[Properties of asymptotic profiles]\label{Thm.Elliptic.Harnack.m}
Assume that $\A$ satisfies (A1), (A2), and (K2). Then there exists a unique positive solution $S$ to the Dirichlet Problem \eqref{Elliptic.prob} with $m>1$. Moreover, let $\sigma$ be as in \eqref{as.sep.var}, and assume that:\\
- either $\sigma=1$ and $2sm\ne \gamma(m-1)$;\\
- or $\sigma<1$ and (K4) holds.\\
Then there exist positive constants $c_0$ and $c_1$ such that the following sharp absolute bounds hold true for all $x\in \Omega$:
\begin{equation}\label{Thm.Elliptic.Harnack.ineq.m.1}
 c_0 \Phi_1(x)^{\sigma/m}\le S(x)\le c_1 \Phi_1(x)^{\sigma/m}\,.
\end{equation}
When $2sm= \gamma(m-1)$ then,  assuming (K4), for all $x\in \Omega$ we have
\begin{equation}\label{Thm.Elliptic.Harnack.ineq.m.1.log}
 c_0 \Phi_1(x)^{1/m}\left(1+|\log\Phi_1(x) |\right)^{1/(m-1)}\le S(x)\le c_1 \Phi_1(x)^{1/m}\left(1+|\log\Phi_1(x) |\right)^{1/(m-1)}\,.\vspace{-1mm}
\end{equation}
\end{thm}
\noindent\textbf{Remark. }As observed in the proof of Theorem \ref{Thm.Asympt},
by applying Theorem \ref{thm.GHP.PME.I} to the separate-variables solution $t^{-\frac{1}{m-1}}S(x)$ we deduce that
\eqref{Thm.Elliptic.Harnack.ineq.m.1} is still true when $\sigma<1$  if, instead of assuming (K4),
we suppose that $K(x,y)\le c_1|x-y|^{-(N+2s)}$ for a.e. $x,y\in \RR^N$
and that $\Phi_1\in C^\gamma(\Omega)$.

\medskip

When $T=0$, the solution $\mathcal U_0$ in \eqref{friendly.giant} is commonly named ``Friendly Giant'', because it takes initial data  $u_0\equiv +\infty$ (in the sense of pointwise limit as $t\to0$) but is bounded for all $t>0$. This term was coined in the study of the standard porous medium equation.

In the following Sections \ref{sect.upperestimates} and \ref{Sec.Lower} we will state and prove our general results concerning upper and lower bounds respectively. These sections are the crux of this paper. The combination of such upper and lower bounds will then be summarized in Section \ref{sect.Harnack}. Consequences of these results in terms of asymptotic behaviour and regularity estimates will be studied in Sections \ref{sec.asymptotic} and \ref{sect.regularity} respectively.\vspace{-2mm}

\section{Upper boundary estimates}\label{sect.upperestimates}
We present a general upper bound that holds under the sole assumptions (A1), (A2), and (K2),
hence valid for all our examples.
\begin{thm}[Absolute boundary estimates]\label{thm.Upper.PME.II}
Let (A1), (A2), and (K2) hold.
Let $u\ge 0$ be a weak dual solution to the (CDP) corresponding to $u_0\in \LL^1_{\Phi_1}(\Omega)$,
and let $\sigma$ be as in \eqref{as.sep.var}. Then, there exists a computable constant $k_1>0$,
depending only on $N, s, m $, and $\Omega$, such that for all $t\ge 0$ and all $x\in \Omega$
\begin{equation}\label{thm.Upper.PME.Boundary.F}
u(t,x) \le \frac{k_1}{t^{\frac1{m-1}}}\,
 \left\{
\begin{array}{ll}
\Phi_1(x_0)^{\sigma/m} &\text{if $\gamma\ne  2sm/(m-1)$},\\
\Phi_1(x_0)^{1/m}\big(1+|\log \Phi_1(x_0)|\big)^{1/(m-1)}&\text{if $\gamma=2sm/(m-1)$}.\\
\end{array}
\right.
\end{equation}
\end{thm}

This absolute bound proves a strong regularization which is independent of the initial datum. It improves the absolute bound in \cite{BV-PPR2-1} in the sense that it exhibits a precise boundary behavior.  The estimate gives the correct behaviour for the solutions $\mathcal{U}_T$ in \eqref{friendly.giant} obtained by separation of variables, see  Theorem \ref{Thm.Elliptic.Harnack.m}. It  turns out that the estimate will be sharp for all nonnegative, nontrivial solutions in the case of the RFL and CFL. We will also see below that the estimate is not always the correct behaviour for the SFL when data are small, as explained in the Introduction
 (see Subsection \ref{ssec.upper.small.data}, and Theorem \ref{prop.counterex} in Section \ref{Sec.Lower}).

\noindent\textit{Proof of Theorem \ref{thm.Upper.PME.II}. }
This subsection is devoted to the proof of Theorems \ref{thm.Upper.PME.II}. The first steps are based on a few basic results of \cite{BV-PPR2-1} that we will also be used in the rest of the paper.

\noindent\textsc{Step 1. Pointwise and absolute upper estimates}\label{sec.upper.partI}

\noindent\textit{Pointwise estimates. } We begin by recalling the basic pointwise estimates which are crucial in the proof of all the upper and lower bounds of this paper.\vspace{-1mm}
\begin{prop}[\cite{BV-PPR1,BV-PPR2-1}]\label{prop.point.est}
It holds\vspace{-1mm}
\begin{equation}\label{thm.NLE.PME.estim.0}
\int_{\Omega}u(t,x)\K(x , x_0)\dx\le \int_{\Omega}u_0(x)\K(x , x_0)\dx \qquad\mbox{for all $t> 0$\,.}
\end{equation}
Moreover, for every $0< t_0\le t_1 \le t$ and almost every $x_0\in \Omega$\,, we have
\begin{equation}\label{thm.NLE.PME.estim}
\frac{t_0^{\frac{m}{m-1}}}{t_1^{\frac{m}{m-1}}}(t_1-t_0)\,u^m(t_0,x_0)
\le \int_{\Omega}\big[u(t_0,x)-u({t_1},x)\big]\K(x , x_0)\dx  \le (m-1)\frac{t^{\frac{m}{m-1}}}{t_0^{\frac{1}{m-1}}}\,u^m(t,x_0)\,.
\end{equation}
\end{prop}
\textit{Absolute upper bounds. } Using the estimates above, in Theorem 5.2 of \cite{BV-PPR2-1}
the authors proved that solutions corresponding to initial data $u_0\in\LL^1_{\Phi_1}(\Omega)$ satisfy\vspace{-1mm}
\begin{equation}\label{thm.Upper.PME.Absolute.F}
\|u(t)\|_{\LL^\infty(\Omega)}\le \frac{K_1}{t^{\frac{1}{m-1}}}\qquad\mbox{for all $t>0$\,,}
\end{equation}
with a constant $K_1$  independent of $u_0$. For this reason, this is called ``absolute bound''.

\noindent\textsc{Step 2. Upper bounds via Green function estimates. }The proof of Theorem \ref{thm.Upper.PME.II} requires the following general statement (see \cite[Proposition 6.5]{BFV-Elliptic}):\vspace{-2mm}
\begin{lem}\label{Lem.Green.2aaa}
Let (A1), (A2), and (K2) hold, and let  $v:\Omega\to \mathbb R$ be a nonnegative bounded function.
Let $\sigma$ be as in \eqref{as.sep.var}, and assume that, for a.e. $x_0\in\Omega$,
\begin{equation}\label{Lem.Green.2.hyp.aaa}
v(x_0)^m\le \kappa_0\int_{\Omega} v(x)\K(x,x_0)\dx.
\end{equation}
Then, there exists a constant $\ka_\infty>0$, depending only on $s,\gamma, m, N,\Omega$, such that the following bound holds true for a.e. $x_0\in \Omega$:
\begin{equation}\label{Lem.Green.2.est.Upper.aaa}
\int_{\Omega} v(x)\K(x,x_0)\dx\le \ka_\infty\kappa_0^{\frac{1}{m-1}}
\left\{
\begin{array}{ll}
\Phi_1(x_0)^\sigma &\text{if $\gamma\ne 2s m/(m-1)$},\\
\Phi_1(x_0)\big(1+|\log \Phi_1(x_0)|^{\frac{m}{m-1}}\big)&\text{if $\gamma= 2s m/(m-1)$}.\\
\end{array}
\right.
\end{equation}\vspace{-2mm}
\end{lem}

\noindent {\sc Step 3. End of the proof of Theorem \ref{thm.Upper.PME.II}. }
We already know that $u(t)\in \LL^\infty(\Omega)$ for all $t>0$ by \eqref{thm.Upper.PME.Absolute.F}.  Also, choosing $t_1=2t_0$ in \eqref{thm.NLE.PME.estim} we deduce that, for $t \ge 0$ and a.e. $x_0\in \Omega$,
\begin{equation}\label{Upper.PME.Step.1.2.c}
u^m(t,x_0) \le \frac{2^{\frac{m}{m-1}}}{t}\int_{\Omega}u(t,x)\K(x , x_0)\dx\,.
\end{equation}
The above inequality corresponds exactly to hypothesis \eqref{Lem.Green.2.hyp.aaa} of Lemma \ref{Lem.Green.2aaa} with the value $\kappa_0=2^{\frac{m}{m-1}}t^{-1}$.  As a consequence, inequality \eqref{Lem.Green.2.est.Upper.aaa} holds,
and we conclude that for a.e. $x_0\in\Omega$ and all $t>0$
\begin{equation}\label{thm.Upper.PME.Boundary.2}
\int_{\Omega}u(t,x)\K(x , x_0)\dx
\le  \frac{\ka_\infty 2^{\frac{m}{(m-1)^2}}}{t^{\frac{1}{m-1}}}\left\{
\begin{array}{ll}
\Phi_1(x_0)^\sigma &\text{if $\gamma\ne \frac{2sm}{m-1}$},\\
\Phi_1(x_0)\big(1+|\log \Phi_1(x_0)|^{\frac{m}{m-1}}\big)&\text{if $\gamma=\frac{2sm}{m-1}$}.\\
\end{array}
\right.
\end{equation}
Hence, combining this bound with \eqref{Upper.PME.Step.1.2.c}, we get
$$
 u^{m}(t,x_0)
\le \frac{k_1^{m}}{t^{\frac{m}{m-1}}}\left\{
\begin{array}{ll}
\Phi_1(x_0)^\sigma &\text{if $\gamma\ne \frac{2sm}{m-1}$},\\
\Phi_1(x_0)\big(1+|\log \Phi_1(x_0)|^{\frac{m}{m-1}}\big)&\text{if $\gamma=\frac{2sm}{m-1}$}.\\
\end{array}
\right.
$$
This proves the upper bounds \eqref{thm.Upper.PME.Boundary.F} and concludes the proof.\qed

\subsection{Upper bounds for small data and small times}\label{ssec.upper.small.data}
As mentioned in the Introduction, the above upper bounds may not be realistic when $\sigma<1$. We have the following estimate for small times if the initial data are sufficiently small.
\begin{thm}\label{thm.Upper.PME.III}
Let $\A$ satisfy (A1), (A2), and (L2). Suppose also that $\A$ has a first eigenfunction $\Phi_1\asymp \dist(x,\partial\Omega)^\gamma$\,, and assume that $\sigma<1$. Finally, we assume that for all $x,y\in \Omega$
\begin{equation}\label{Operator.Hyp.upper.III}\begin{split}
K(x,y)\le \frac{c_1}{|x-y|^{N+2s}}
\left(\frac{\Phi_1(x)}{|x-y|^\gamma }\wedge 1\right)
\left(\frac{\Phi_1(y)}{|x-y|^\gamma }\wedge 1\right)
&\quad\mbox{ and }\quad B(x)\le c_1\Phi_1(x)^{-\frac{2s}{\gamma}}\,.
\end{split}\end{equation}
Let $u\ge 0$ be a weak dual solution to (CDP) corresponding to $u_0\in\LL^1_{\Phi_1}(\Omega)$.
Then, for every initial data $u_0\leq A\,\Phi_1^{1-2s/\gamma}$ for some $A>0$, we have
$$
u(t)\leq \frac{\Phi_1^{1-\frac{2s}{\gamma}}}{[A^{1-m} -\tilde C t]^{m-1}}
\qquad \text{on }[0,T_A],\qquad\text{where } T_A:=\frac{1}{\tilde CA^{m-1}},
$$
and the constant $\tilde C>0$, that depends only on $N,s,m, \lambda_1, c_1$, and $\Omega$.
\end{thm}

\noindent\textbf{Remark. }This result applies to the SFL. Notice that when $\sigma<1$ we have always $1-\frac{2s}{\gamma}>\sigma/m$\,, hence in this situation small data have a smaller behaviour at the boundary than the one predicted in Theorem \ref{thm.Upper.PME.II}. This is not true for ``big'' data, for instance for solution  obtained by separation of variables, as already said.

\noindent\textit{Proof of Theorem \ref{thm.Upper.PME.III}. }
In view of our assumption on the initial datum, namely $u_0\le A\,\Phi_1^{1-2s/\gamma}$, by comparison it is enough to prove that
 the function
\[
\ua(t,x)=F(t)\Phi_1(x)^{1-\frac{2s}{\gamma}},\qquad F(t)=\frac{1}{[A^{1-m} -\tilde Ct]^{m-1}}
\]
is a supersolution (i.e., $\partial_t\ua\ge -\A\ua^m$) in $(0,T_A)\times \Omega$
provided we choose $\tilde C$ sufficiently large.

To this aim, we use the following elementary inequality, whose proof is left to the interested reader:
for any $\eta>1$ and any $M>0$ there exists $\widetilde{b}=\widetilde{b}(M)>0$   such that letting $\widetilde{\eta}:=\eta\wedge 2$
\begin{equation}\label{supersol.1.step3}
a^\eta-b^\eta\le \eta\,b^{\eta-1}(a-b)+  \widetilde{b} |a-b|^{\widetilde{\eta}}\,,\qquad\mbox{ for all $0\le a,b\le M$.}
\end{equation}
We apply inequality \eqref{supersol.1.step3} to $a=\Phi_1(y)$ and $b=\Phi_1(x)$, $\eta=m(1-\frac{2s}{\gamma})$, noticing that $\eta>1$ if and only if $\sigma<1$\,, and we obtain (recall that $\Phi_1$ is bounded)
\begin{equation*}
\begin{split}
 \ua^m(t,y)-\ua^m(t,x)&= F(t)^m \left(\Phi_1(y)^{m(1-\frac{2s}{\gamma})}-\Phi_1(x)^{m(1-\frac{2s}{\gamma})} \right)
 = F(t)^m \left(\Phi_1(y)^\eta-\Phi_1(x)^\eta \right)\\
&\le \eta\,F(t)^m\Phi_1(x)^{\eta-1}\left[\Phi_1(y)-\Phi_1(x)\right]
 +\widetilde{b}\, F(t)^m\left|\Phi_1(y)-\Phi_1(x)\right|^{\widetilde{\eta}}\\
&\le \eta\,F(t)^m\Phi_1(x)^{\eta-1} \left[\Phi_1(y)-\Phi_1(x)\right]
 + \widetilde{b}\,F(t)^m c_\gamma^{\widetilde{\eta}} |x-y|^{\widetilde{\eta}\gamma},
\end{split}
\end{equation*}
where in the last step we have used that $\left|\Phi_1(y)-\Phi_1(x)\right|\le c_\gamma |x-y|^\gamma$.
Since $B\le c_1\Phi_1^{-2s/\gamma}$,

\[\begin{split}
\int_{\RR^N}\left[\Phi_1(y)-\Phi_1(x)\right]K(x,y)\dy
&=-\A \Phi_1(x) + B(x)\Phi_1(x)
\le -\lambda_1\Phi_1(x) + c_1\Phi_1(x)^{1-\frac{2s}{\gamma}}\,, \\
\end{split}\]
Thus, recalling that $\eta,\widetilde{\eta}>1$ and that $\Phi_1$ is bounded, it follows
\begin{equation}\label{supersol.2}\begin{split}
-\A[\ua^m](x)   &=\int_{\RR^N}\left[ \ua^m(t,y)-\ua^m(t,x)\right]K(x,y)\dy + B(x)\ua^m(t,x)\\
&\le \eta\,F(t)^m\Phi_1(x)^{\eta-1} \left[-\lambda_1\Phi_1(x) + c_1\Phi_1(x)^{1-\frac{2s}{\gamma}}\right]  + B(x)F(t)^m \Phi_1^{\eta}(x)\\
&+ \widetilde{b}\,c_\gamma^{\widetilde{m}}\,F(t)^m  \int_{\RR^N}|x-y|^{\widetilde{\eta}\gamma}K(x,y)\dy\\
&\le \widetilde{c}F(t)^m \left(\Phi_1(x)^{\eta-\frac{2s}{\gamma}} + \int_{\RR^N}|x-y|^{\widetilde{\eta}\gamma}K(x,y)\dy\right)
\end{split}
\end{equation}
Next, we claim that, as a consequence of \eqref{Operator.Hyp.upper.III})
\begin{equation}\label{supersol.3}
\int_{\RR^N}|x-y|^{\widetilde{\eta}\gamma}K(x,y)\dy
\le c_4\Phi_1(x)^{1-\frac{2s}{\gamma}}\,.
\end{equation}
Postponing for the moment the proof of the above inequality, we first show how conclude: combining \eqref{supersol.2} and \eqref{supersol.3} we have
\[
-\A\ua^m \le c_5 F(t)^m \Phi_1(x)^{1-\frac{2s}{\gamma}}= F'(t)\Phi_1(x)^{1-\frac{2s}{\gamma}}= \partial_t \ua
\]
where we used that
$F'(t)=c_5 F(t)^{\widetilde{m}}$
provided $\Tilde C=c_5(m-1)$. This proves that $\ua$ is a supersolution in $(0,T)\times \Omega$.
Hence the proof is concluded once we prove inequality \eqref{supersol.3}; for this, using hypothesis \eqref{Operator.Hyp.upper.III} and choosing $r=\Phi_1(x)^{1/\gamma}$ we have
\[\begin{split}
\int_{\RR^N}&|x-y|^{\widetilde{\eta}\gamma}K(x,y)\dy
\le c_1\int_{B_r(x)}\frac{1}{|x-y|^{N+2s-\widetilde{\eta}\gamma}}\dy+c_1\Phi_1(x)\int_{\Omega\setminus B_r(x)}\frac{1}{|x-y|^{N+2s+\gamma-\widetilde{\eta}\gamma}} \dy\\
&\le c_2 r^{\widetilde{\eta}\gamma-2s}+c_1\frac{\Phi_1(x)}{r^{2s}}
    \int_{\Omega\setminus B_r(x)}\frac{1}{|x-y|^{N+\gamma-\widetilde{\eta}\gamma}}\dy
    = c_2 r^{\widetilde{\eta}\gamma-2s}+ c_3\frac{\Phi_1(x)}{r^{2s}} \le c_4\Phi_1(x)^{1-\frac{2s}{\gamma}}\,,
\end{split}\]
where we used that $\widetilde{\eta}\gamma-2s>0$ and $\tilde\eta>1$.\qed

\noindent{\bf Remark.} For operators for which the previous assumptions hold with  $B\equiv 0$, we can actually prove a better upper bound for ``smaller data'', namely:
\begin{cor}\label{thm.Upper.PME.IV}Under the assumptions of Theorem \ref{thm.Upper.PME.III}, assume that moreover $B\equiv 0$ and  $u_0\leq A\,\Phi_1$ for some $A>0$. Then, we have
$$
u(t)\leq \frac{\Phi_1}{[A^{1-m} -\tilde C t]^{m-1}}
\qquad \text{on }[0,T_A],\qquad\text{where } T_A:=\frac{1}{\tilde CA^{m-1}},
$$
and the constant $\tilde C>0$, that depends only on $N,s,m, \lambda_1, c_1$, and $\Omega$.
\end{cor}
\noindent {\bf Proof.~}We have to show that $\ua(t,x)=F(t)\Phi_1(x)$ is a supersolution: we essentially repeat the proof of Theorem \ref{thm.Upper.PME.III} with $\gamma=m$ (formally replace $1-2s/\gamma$ by $1$), taking into account that $B\equiv 0$ and $u_0\leq A\,\Phi_1$.\qed


\section{Lower bounds}\label{Sec.Lower}

This section is devoted to the proof of all the lower bounds summarized later in the  main Theorems \ref{thm.GHP.PME.I}, \ref{thm.GHP.PME.II}, and \ref{thm.GHP.PME.III}. The general situation is quite involved to describe, so we will separate several cases and we will indicate for which examples it holds for the sake of clarity.

\medskip

\noindent $\bullet$ \textbf{Infinite speed of propagation: universal lower bounds. }First, we are going to quantitatively establish that all nonnegative weak dual solutions of our problems are in fact positive in $\Omega$ for all $t>0$. This result is valid for all nonlocal operators considered in this paper.

\begin{thm}\label{thm.Lower.PME}
Let $\A$ satisfy (A1), (A2), and (L2). Let $u\ge 0$ be a weak dual solution to the (CDP) corresponding to $u_0\in \LL^1_{\Phi_1}(\Omega)$. Then there exists a constant $\kb_0>0$ such that the following inequality holds:
\begin{equation}\label{thm.Lower.PME.Boundary.1}
u(t,x)\ge \kb_0\,\left(1\wedge \frac{t}{t_*}\right)^{\frac{m}{m-1}}\frac{\Phi_1(x)}{t^{\frac{1}{m-1}}}\qquad\mbox{for all $t>0$ and a.e. $x\in \Omega$}\,.
\end{equation}
Here $t_*=\k_*\|u_0\|_{\LL^1_{\Phi_1}(\Omega)}^{-(m-1)}$, and the constants $\kb_0$ and $\k_*$ depend only on   $N,s,\gamma, m, c_0,c_1$, and $\Omega$\,.
\end{thm}
Notice that, for $t
\geq t_*$, the dependence on the initial data disappears from the lower bound,
as inequality reads
$$
u(t)\geq \kb_0 \frac{\Phi_1}{t^{\frac{1}{m-1}}}\qquad \forall\,t \geq t_*,
$$
where $\kb_0$ is an absolute constant.
Assumption (L2) on the kernel $K$ of $\A$ holds for all examples mentioned in Section \ref{sec.examples}.

Clearly, the power in this lower bound does not match the one of the general upper bounds of Theorem \ref{thm.Upper.PME.II}, hence we can not expect these bounds to be sharp. However, when $\sigma<1$, for small times and small data and when $B\equiv 0$, the lower bounds \eqref{thm.Lower.PME.Boundary.1} match the upper bounds of Corollary \ref{thm.Upper.PME.IV}, hence they are sharp. Theorem \ref{thm.Lower.PME} shows that, even in the ``worst case scenario'', there is a quantitative lower bound for all positive times, and shows infinite speed of propagation.

\noindent$\bullet$ \textbf{Matching lower bounds I. }Actually, in many cases the kernel of the nonlocal operator satisfies a stronger property, namely $\inf_{x,y\in \Omega}K(x,y)\ge \underline{\kappa}_\Omega>0$  and $B\equiv 0$, in which case we can actually obtain sharp lower bounds for all times. Here we do not consider the potential logarithmic correction that may appear in ``critical case'' $2sm= \gamma(m-1)$: indeed, as far as examples are concerned, the next Theorem applies to the RFL and the CFL, for which   $2sm> \gamma(m-1)$.
\begin{thm}\label{Thm.lower.B}
Let $\A$ satisfy (A1), (A2), and (L1).
Furthermore, suppose that $\A$ has a first eigenfunction $\Phi_1\asymp \dist(x,\partial\Omega)^\gamma$\,. Let $\sigma$ be as in \eqref{as.sep.var} and assume that:\\
- either $\sigma=1$;\\
- or $\sigma<1$, $K(x,y)\le c_1|x-y|^{-(N+2s)}$ for a.e. $x,y\in \RR^N$,
and $\Phi_1\in C^\gamma(\overline{\Omega})$.\\
Let $u\ge 0$ be a weak dual solution to the (CDP) corresponding to $u_0\in \LL^1_{\Phi_1}(\Omega)$. Then there exists a constant $\kb_1>0$ such that the following inequality holds:
\begin{equation}\label{Thm.B.lower.bdd}
u(t,x)\ge \kb_1 \left(1\wedge \frac{t}{t_*}\right)^{\frac{m}{m-1}}\frac{\Phi_1(x)^{\sigma/m}}{t^{\frac{1}{m-1}}} \qquad\mbox{for all $t>0$ and a.e. $x\in \Omega$}\,,
\end{equation}
where $t_*=\k_*\|u_0\|_{\LL^1_{\Phi_1}(\Omega)}^{-(m-1)}$.  The constants $\k_*$ and $\kb_1$ depend only on   $N,s,\gamma, m, \kb_\Omega,c_1,\Omega,$ and $\|\Phi_1\|_{C^\gamma(\Omega)}$.
\end{thm}

\noindent\textbf{Remarks. }(i) As in the case of the Theorem \ref{thm.Lower.PME}, for large times the dependence on the initial data disappears from the lower bound and we have absolute lower bounds.

\noindent(ii) The boundary behavior is sharp  when $2sm\ne \gamma(m-1)$  in view of the upper bound from Theorem \ref{thm.Upper.PME.II}.

\noindent(iii) This theorem applies to the RFL and the CFL, but not to the SFL (or, more in general, spectral powers of elliptic operators), see Sections \ref{ssec.examples} and \ref{sec.hyp.L}. In the case of the RFL, this result was obtained in Theorem 1 of \cite{BFR}.

\medskip

We have already seen the example of the separate-variables solutions \eqref{friendly.giant}
that have a very definite behavior at the boundary $\partial \Omega$. The analysis of general solutions leads to completely different situations for
$\sigma=1$ and $\sigma<1$.

\medskip

\noindent $\bullet$ \textbf{Matching lower bounds II. The case $\sigma=1$. }When $\sigma=1$ we can establish  a quantitative lower bound near the boundary that matches the separate-variables behavior for large times (except in the case $2sm= \gamma(m-1)$ where the result is false, see Theorem \ref{prop.counterex} below). We do not need the assumption of non-degenerate kernel, so SFL can be considered.

\begin{thm}\label{thm.Lower.PME.large.t}
Let $(A1)$, $(A2)$, and $(K2)$ hold, and let $\sigma=1$. Let $u\ge 0$ be a weak dual solution to the (CDP) corresponding to $u_0\in \LL^1_{\Phi_1}(\Omega)$. There exists a constant $\kb_2>0$ such that
\begin{equation}\label{thm.Lower.PME.Boundary.large.t}
 u(t,x) \ge \kb_2\,\frac{\Phi_1(x)^{1/m}}{t^{\frac{1}{m-1}}}\qquad\mbox{for all $t\ge t_*$ and a.e. $x\in \Omega$}\,.
\end{equation}
Here, $t_*=\k_*\|u_0\|_{\LL^1_{\Phi_1}(\Omega)}^{-(m-1)}$, and the constants $\k_* $ and $\kb_2$ depend only on   $N,s,\gamma, m $, and $\Omega$\,.
\end{thm}
\noindent {\bf Remarks.} (i) At first sight, this theorem may seem weaker than the previous positivity result. However, this result  has wider applicability since it holds  under the only assumption (K2)
on $\K.$  In particular it is valid in the local case $s=1$, where the finite speed of propagation makes it impossible to have global lower bounds for small times.

\noindent (ii) When $\A=-\Delta$ the result has been proven in \cite{Ar-Pe} and \cite{JLVmonats}
by quite different methods. On the other hand, our method is very general and immediately applies to the case when $\A$ is an elliptic operator with $C^1$ coefficients, see Section \ref{sec.examples}.

\noindent  (iii) This result fixes a small error in Theorem 7.1 of \cite{BV-PPR1} where the power $\sigma$ was not present.  

\medskip

\noindent$\bullet$ \textbf{The anomalous lower bounds with small data. }As shown in Theorem \ref{thm.Lower.PME},
the lower bound $u(t)\gtrsim \Phi_1$ is always valid. We now discuss the possibility of improving this bound.

Let $S$  solve the elliptic problem \eqref{Elliptic.prob}.
It follows by comparison whenever $u_0
\geq \epsilon_0 S$ with $\epsilon_0>0$ then $u(t)\geq \frac{S}{(T_0+t)^{1/(m-1)}}$, where $T_0=\epsilon_0^{1-m}$.
Since $S\asymp \Phi_1^{\sigma/m}$ under (K4) (up to a possible logarithmic correction in the critical case, see Theorem \ref{Thm.Elliptic.Harnack.m}), there are initial data for which the lower behavior is dictated by $\Phi_1(x)^{\sigma/m}t^{-1/(m-1)}$.
More in general,
as we shall see in Theorem \ref{Thm.Asympt.0}, given any initial datum $u_0
\in \LL^1_{\Phi_1}(\Omega)$ the function $v(t,x):=t^{\frac{1}{m-1}}u(t)$ always converges to $S$ in $\LL^\infty(\Omega)$ as $t\to \infty$, independently of the value of $\sigma$.
Hence, one may conjecture that there should exist a waiting time $t_*>0$ after which the lower behavior is dictated by $\Phi_1(x)^{\sigma/m}t^{-1/(m-1)}$, in analogy with what happens for the classical porous medium equation.
As we shall see, this is actually {\em false} when $\sigma<1$ or $2sm= \gamma(m-1)$. 
Since for large times $v(t,x)$ must look like $S(x)$ in uniform norm away from the boundary (by the interior regularity that we will prove later), the contrasting situation for large times could be described as `dolphin's head' with the `snout' flatter than the `forehead'. As $t\to\infty$ the forehead progressively fills the whole domain.

The next result shows that, in general, we cannot hope to prove that $u(t)$
is larger than $\Phi_1^{1/m}$. In particular, when $\sigma<1$ or $2sm= \gamma(m-1)$, this shows that the behavior $u(t)\asymp S$ cannot hold.
\begin{thm}\label{prop.counterex}
Let (A1), (A2), and (K2) hold, and $u\ge 0$ be a weak dual solution to the (CDP) corresponding to a nonnegative initial datum $u_0\in \LL^1_{\Phi_1}(\Omega)$.
Assume that $u_0(x)\le C_0\Phi_1(x)$ a.e. in $\Omega$ for some $C_0>0$. Then
there exists a constant $\hat\kappa$, depending only $N,s,\gamma, m $, and $\Omega$, such that
$$
u(t,x)^m \leq C_0\hat\kappa \frac{\Phi_1(x)}{t}\qquad\mbox{for all $t>0$ and a.e.
$x\in \Omega$}\,.
$$
In particular, if $\sigma<1$ (resp. $2sm= \gamma(m-1)$), then
$$
\lim_{x\to \partial\Omega}\frac{u(t,x)}{\Phi_1(x)^{\sigma/m}}= 0 \quad \Bigl(\text{resp.} \lim_{x\to \partial\Omega}\frac{u(t,x)}{\Phi_1(x)^{1/m}\left(1+|\log\Phi_1(x) |\right)^{1/(m-1)}}= 0\Bigr)\qquad \text{for any $t>0$.}
$$
\end{thm}

The proposition above could make one wonder whether the sharp general lower bound
could be given by $\Phi_1^{1/m}$, as in the case $\sigma=1$.
Recall that, under rather minimal assumptions on the kernel $K$
associated to $\A$, we have a universal lower bound for $u(t)$ in terms of $\Phi_1$
(see Theorem \ref{thm.Lower.PME}).
Here we shall see that, under (K4), the bound $u(t)\gtrsim \Phi_1^{1/m}$
is false for $\sigma<1$.

\begin{thm}\label{prop.counterex2}
Let (A1), (A2), and (K4) hold, and let $u\ge 0$ be a weak dual solution to the (CDP) corresponding to a nonnegative initial datum $u_0 \leq C_0 \Phi_1$ for some $C_0>0$.
Assume that there exist constants $\kb,T,\alpha>0$ such that
$$
u(T,x)\geq \kb\Phi_1^\alpha(x)\qquad\mbox{for a.e. $x\in \Omega$}\,.
$$
Then $\alpha \geq 1-\frac{2s}{\gamma}$.
In particular $\alpha>\frac{1}{m}$ if $\sigma<1$.
\end{thm}
We devote the rest of this section to the proof of the above results, and to this end we collect in the first two subsections some preliminary lower bounds and results about approximate solutions.
\subsection{Lower bounds for weighted norms}\label{Sect.Weighted.L1}
Here we prove some useful lower bounds for weighted norms, which follow from the $\LL^1$-continuity for ordered solutions in the version proved in Proposition 8.1 of \cite{BV-PPR2-1}.
\begin{lem}[Backward in time $\LL^1_{\Phi_1}$ lower bounds]\label{cor.abs.L1.phi} Let $u$ be a solution to (CDP) corresponding to the initial datum $u_0\in\LL^1_{\Phi_1}(\Omega)$. For all
\begin{equation}\label{L1weight.PME.estimates.2}
0\le \tau_0\le  t\le \tau_0+ \frac{1}{\big(2\bar K\big)^{1/(2s\vartheta_{\gamma})}\|u(\tau_0)\|_{\LL^1_{\Phi_1}(\Omega)}^{m-1}}
\end{equation}
we have
\begin{equation}\label{L1weight.PME.estimates.3}
\frac{1}{2}\int_{\Omega}u(\tau_0,x)\Phi_1(x)\dx\le \int_{\Omega}u(t,x)\Phi_1(x)\dx\,,
\end{equation}
where $\vartheta_{\gamma}:=1/[2s+(N+\gamma)(m-1)]$ and $\bar K>0$ is a computable constant.
\end{lem}
\noindent\textsl{Proof of Lemma \ref{cor.abs.L1.phi}. }We recall the inequality of Proposition 8.1 of \cite{BV-PPR2-1}, adapted to our case: for all $0\le \tau_0\le \tau,t$ we have
\begin{equation}\label{L1weight.contr.estimates.1}
\int_{\Omega} u(\tau,x) \Phi_1(x)\dx
\le \int_{\Omega} u(t,x) \Phi_1(x)\dx
+ \bar K\|u(\tau_0)\|_{\LL^1_{\Phi_1}(\Omega)}^{2s(m-1) \vartheta_{\gamma}+1}\,\left|t-\tau\right|^{2s\vartheta_{\gamma}} \,.
\end{equation}
Choosing $\tau=\tau_0$ in the above inequality, we get
\begin{equation}\label{cor.0.1}\begin{split}
\left[1- K_9\|u(\tau_0)\|_{\LL^1_{\Phi_1}(\Omega)}^{2s(m-1) \vartheta_{\gamma}}\,\left|t-\tau_0\right|^{2s\vartheta_{\gamma}}\right]\int_{\Omega}u(\tau_0,x)\Phi_1(x)\dx &\le \int_{\Omega}u(t,x)\Phi_1(x)\dx\,.
\end{split}
\end{equation}
Then \eqref{L1weight.PME.estimates.3} follows from \eqref{L1weight.PME.estimates.2}\,.\qed

We also need a lower bound for $\LL^p_{\Phi_1}(\Omega)$ norms.
\begin{lem}\label{lem.abs.lower}
Let $u$ be a solution to (CDP) corresponding to the initial datum $u_0\in\LL^1_{\Phi_1}(\Omega)$. Then the following lower bound holds true for any $t\in [0,t_*]$ and $p\ge 1$:
\begin{equation}\label{lem1.lower.bdd}
c_2 \left(\int_{\Omega}u_0(x)\Phi_1(x)\dx \right)^p \le \int_{\Omega}u^p(t,x)\Phi_1(x)\dx 
\end{equation}
Here $t_*=c_*\|u_0\|_{\LL^1_{\Phi_1}(\Omega)}^{-(m-1)}$, where $c_2, c_*>0$ are positive constants that depend only on $N,s,m,p,\Omega$.
\end{lem}
The proof of this Lemma is an easy adaptation of the proof of Lemma 2.2 of \cite{BFR}\,, so we skip it.  Notice that $c_*$ has explicit form given in {\rm \cite{BV-PPR1,BV-PPR2-1, BFR}}, while the form of $c_2$ is given in the proof of Lemma 2.2 of \cite{BFR}.
 \subsection{Approximate solutions}
 \label{sect:approx sol}
To prove our lower bounds, we will need a special class of approximate solutions $u_\delta$. We will list now the necessary details. In the case when $\A$ is the Restricted Fractional Laplacian (RFL) (see Section \ref{sec.examples}) these solutions have been used in the Appendix II of \cite{BFR}, where complete proofs can be found; the proof there holds also for the operators considered here. The interested reader can easily adapt the proofs in \cite{BFR} to the current case.

Let us fix $\delta>0$ and consider the problem:
\begin{equation}\label{probl.approx.soln}
\left\{
\begin{array}{lll}
\partial_t v_\delta=-\A\left[(v_\delta+\delta)^m -\delta^m \right]& \qquad\mbox{for any $(t,x)\in (0,\infty)\times \Omega$}\\
v_\delta(t,x)=0 &\qquad\mbox{for any $(t,x)\in (0,\infty)\times (\RR^N\setminus\Omega)$}\\
v_\delta(0,x)=u_0(x) &\qquad\mbox{for any $x\in\Omega$}\,.
\end{array}
\right.
\end{equation}
Next, we define
\[
u_\delta:=v_\delta+\delta.
\]
We summarize here below the basic properties of $\ud$.

Approximate solutions $\ud$ exist, are unique, and bounded for all $(t,x)\in (0,\infty)\times\overline{\Omega}$ whenever $0\le u_0\in \LL^1_{\Phi_1}(\Omega)$\,.
Also, they are uniformly positive: for any $t \geq 0$,
\begin{equation}\label{approx.soln.positivity}
\ud(t,x)\ge \delta>0 \qquad\mbox{for a.e. $x \in \Omega$.}
\end{equation}
This implies that the equation for $\ud$ is never degenerate in the interior,
so solutions are smooth as the linear parabolic theory with the kernel $K$
allows them to be (in particular, in the case of the fractional laplacian, they are $C^\infty$ in space and $C^1$ in time).
Also, by a comparison principle,
for all $\delta>\delta'>0$ and $t \geq 0$,
\begin{equation}\label{approx.soln.comparison1}
\ud(t,x)\ge u_{\delta'}(t,x) \qquad\mbox{for $x \in \Omega$}
\end{equation}
and
\begin{equation}\label{approx.soln.comparison}
\ud(t,x)\ge u(t,x) \qquad\mbox{for a.e. $x \in \Omega$\,.}
\end{equation}
Furthermore, they converge in $\LL^1_{\Phi_1}(\Omega)$ to $u$ as $\delta \to 0$:
\begin{equation}\label{Lem.0}
\|\ud(t)-u(t)\|_{\LL^1_{\Phi_1}(\Omega)}\le \|\ud(0)-u_0\|_{\LL^1_{\Phi_1}(\Omega)}=\delta\,\|\Phi_1\|_{\LL^1(\Omega)}\,.
\end{equation}
As a consequence of \eqref{approx.soln.comparison1} and \eqref{Lem.0}, we deduce that $\ud$ converge pointwise to $u$
at almost every point: more precisely, for all $t\geq 0,$
\begin{equation}\label{limit.sol}
{u}(t,x)=\lim_{\delta\to 0^+}\ud(t,x)\qquad\mbox{for a.e. $x \in \Omega$\,.}
\end{equation}

\subsection{Proof of Theorem \ref{thm.Lower.PME}}The proof consists in showing that
\[
u(t,x) \geq \ub(t,x):=k_0\, t \,\Phi_1(x)
\]
for all $t\in [0,t_*]$, where the parameter $k_0>0$ will be fixed later.
Note that, once the inequality $u
\geq \ub$ on $[0,t_*]$ is proved, we conclude as follows: since
$t\mapsto t^{\frac{1}{m-1}}\,u(t,x)$ is nondecreasing in $t>0$ for a.e. $x\in \Omega$
(cf. (2.3) in \cite{BV-PPR2-1}) we have
$$u(t,x)\ge \left( \frac{t_*}{t}\right)^{\frac{1}{m-1}}u(t_*,x) \geq k_0\, t_*  \left( \frac{t_*}{t}\right)^{\frac{1}{m-1}}\,\Phi_1(x)\qquad \text{ for all $t\ge t_*$\,.}$$
Then, the result will follow $\kb_0=k_0t_*^{\frac{m}{m-1}}$
(note that, as we shall see below, $k_0t_*^{\frac{m}{m-1}}$ can be chosen independently of $u_0$).
Hence, we are left with proving that $u
\geq \ub$ on $[0,t_*]$.

\noindent$\bullet~$\textsc{Step 1. }\textit{Reduction to an approximate problem. }Let us fix $\delta>0$ and consider the approximate solutions $\ud$ constructed in Section \ref{sect:approx sol}.
We shall prove that  $\ud
\geq \ub$ on $[0,t_*]$, so that the result will follow by the arbitrariness of $\delta$.

\noindent$\bullet~$\textsc{Step 2. }\textit{We claim that $\ub(t,x)< \ud(t,x)$ for all $0\le t\le t_*$ and $x\in \Omega$, for a suitable choice of $k_0>0$\,.} Assume that the inequality $\ub<\ud$ is false in $[0,t_*]\times \overline\Omega$, and let $(t_c,x_c)$ be the first contact point between $\ub$ and $\ud$. Since $\ud=\delta>0=\ub$ on the lateral boundary, $(t_c,x_c)\in(0,t_*]\times\Omega$\,. Now, since $(t_c,x_c)\in (0,t_*]\times\Omega$ is the first contact point, we necessarily have that
\begin{equation}\label{contact.1}
\ud(t_c,x_c)= \ub(t_c,x_c)\qquad\mbox{and}\qquad
\ud(t,x)\ge \ub(t,x)\,\quad\forall t\in [0,t_c]\,,\;\; \forall x\in\overline{\Omega}\,.
\end{equation}
Thus, as a consequence,
\begin{equation}\label{contact.2}
\partial_t \ud(t_c,x_c)\le \partial_t \ub(t_c,x_c)=k_0\,\Phi_1(x_c)\,.
\end{equation}
Next, we observe the following  Kato-type inequality holds: for any nonnegative function $f$,
\begin{equation}\label{Kato.ineq}
\A(f^m)\le mf^{m-1}\A f.
\end{equation}
Indeed, by convexity, $f(x)^m-f(y)^m \le m [f(x)]^{m-1}(f(x)-f(y))$, therefore
\begin{equation*}
\begin{split}
\A(f^m)(x)&=\int_{\RR^N}[f(x)^m-f(y)^m]\,K(x,y)\dy  + B(x)f(x)^m \\
&\le m [f(x)]^{m-1}\int_{\RR^N}[f(x)-f(y)]\,K(x,y)\dy  + B(x)f(x)^m \\
&= m [f(x)]^{m-1}\left[\int_{\RR^N}[f(x)-f(y)]\,K(x,y)\dy + B(x)f(x) \right] - (m-1)B(x)f(x)^m\\
&\le m [f(x)]^{m-1}\A f(x)\,.
\end{split}
\end{equation*}
As a consequence of \eqref{Kato.ineq}, since $t_c\le t_*$ and $\Phi_1$ is bounded,
\begin{multline}\label{contact.2b}
\A(\ub^m)(t,x)\le m\ub^{m-1}\A(\ub)=m [k_0 t \Phi_1(x)]^{m-1}\,k_0 t \A(\Phi_1)(x)\\
=m \lambda_1[k_0 t \Phi_1(x)]^m \le \k_1 (t_*k_0)^m\Phi_1(x)\,,
\end{multline}
Then, using  \eqref{contact.2} and \eqref{contact.2b}, we establish an upper bound for $-\A(u_\delta^m-\ub^m)(t_c,x_c)$ as follows:
\begin{equation}\label{contact.3}
-\A[u_\delta^m- \ub^m](t_c,x_c)=\partial_t \ud(t_c,x_c)+\A(\ub^m)(t_c,x_c)
\le k_0\, \left[1+\k_1t_*^mk_0^{m-1}\right]\Phi_1(x_c).
\end{equation}
Next, we want to prove lower bounds for $-\A(u_\delta^m-\psi^m)(t_c,x_c)$, and this is the point where the nonlocality of the operator enters, since we make essential use of hypothesis (L2). We recall that by  \eqref{contact.1} we have $u_\delta^m(t_c,x_c)=\ub^m(t_c,x_c)$, so that assumption (L2) gives
\begin{equation*}\begin{split}
-\A &\left[u_\delta^m-\ub^m\right](t_c,x_c)= -\A \left[u_\delta^m-\ub^m\right](t_c,x_c) + B(x_c)[u_\delta^m(t_c,x_c)-\ub^m(t_c,x_c)]\\
&=-\int_{\RR^N}\left[\big(u_\delta^m(t_c,x_c)-u_\delta^m(t_c,y)\big)-\big(\ub^m(t_c,x_c)-\ub^m(t_c,y)\big)\right]K(x_c,y)\dy\\
&=\int_{\Omega}\left[u_\delta^m(t_c,y)-\ub^m(t_c,y)\right]K(x_c,y)\dy
\ge c_0\Phi_1(x_c)\int_{\Omega}\left[u_\delta^m(t_c,y)-\ub^m(t_c,y)\right]\Phi_1(y)\dy\,,\\
\end{split}
\end{equation*}
from which it follows (since $\ub^m= [k_0 t \Phi_1(x)]^m   \le \k_2(t_*k_0)^m $)
\begin{equation}\label{contact.44}\begin{split}
-\A &\left[u_\delta^m-\ub^m\right](t_c,x_c)\\
&\ge c_0\Phi_1(x_c)\int_{\Omega}u_\delta^m(t_c,y)\Phi_1(y)\dy-c_0\Phi_1(x_c)\int_{\Omega}\ub^m(t_c,y)\Phi_1(y)\dy.\\
&\ge c_0\Phi_1(x_c)\int_{\Omega}u_\delta^m(t_c,y)\Phi_1(y)\dy-c_0\Phi_1(x_c)\k_3\,(t_*k_0)^m.
\end{split}
\end{equation}
Combining the upper and lower bounds \eqref{contact.3} and \eqref{contact.44} we obtain
\begin{equation}\label{contact.5}\begin{split}
c_0\Phi_1(x_c)\int_{\Omega}u_\delta^m(t_c,y)\Phi_1(y)\dy &\le k_0\, \left[1+ (\k_1+\k_3)t_*^m k_0^{m-1}  \right]\Phi_1(x_c)\,.
\end{split}
\end{equation}
Hence, recalling \eqref{lem1.lower.bdd}, we get
$$
c_2 \left(\int_{\Omega}u_0(x)\Phi_1(x)\dx \right)^m\le  \int_{\Omega}u_\delta^m(t_c,y)\Phi_1(y)\dy
\le \frac{k_0}{c_0}\, \left[1+ (\k_1+\k_3)t_*^mk_0^{m-1}  \right].
$$
Since $t_*=\k_*\|u_0\|_{\LL^1_{\Phi_1}(\Omega)}^{-(m-1)}$, this yields
$$
c_2\k_*^{\frac{m}{m-1}}t_*^{-\frac{m}{m-1}} \le \frac{k_0}{c_0}\, \left[1+ (\k_1+\k_3)t_*^mk_0^{m-1}  \right]
$$
which gives the desired contradiction provided we choose $k_0$ so that $\kb_0:=k_0t_*^{\frac{m}{m-1}}$ is universally small.
\qed

\subsection{Proof of Theorem \ref{Thm.lower.B}. }The proof proceeds along the lines of the proof of Theorem \ref{thm.Lower.PME}, so we will just briefly mention the common parts.

We want to show that
\begin{equation}\label{lower.barrier}
\ub(t,x):=\kappa_0\,t\, \Phi_1(x)^{\sigma/m}\,,
\end{equation}
is a lower barrier for our problem on $[0,t_*]\times \Omega$ provided $\kappa_0$ is small enough.
More precisely, as in the proof of Theorem \ref{thm.Lower.PME}, we aim to prove that $\ub<\ud$ on $[0,t_*]$, as the lower bound for $t \geq t_*$ then follows by monotonicity.

Assume by contradiction that the inequality $\ub(t,x)< u_{\delta}(t,x)$
is false inside $[0,t_*]\times \overline\Omega$.
Since $\ub<\ud$ on the parabolic boundary, letting $(t_c,x_c)$ be the first contact point, we necessarily have that $(t_c,x_c)\in (0,t_*]\times\Omega$. The desired contradiction will be obtained by combining the upper and lower bounds (that we prove below) for the quantity $-\A\left[u_\delta^m- \ub^m\right](t_c,x_c)$ , and then choosing $\kappa_0>0$ suitably small. In this direction, it is convenient in what follows to assume that
\begin{equation}\label{k0.first.cond}
\k_0\le 1\wedge t_*^{-\frac{m}{m-1}}\qquad\mbox{so that}\qquad \k_0^{m-1}t_*^m\le 1\,.
\end{equation}

\noindent\textit{Upper bound. }We first establish the following upper bound: there exists a constant $\overline{A}>0$ such that
\begin{equation}\label{contact.up}
-\A\left[u_\delta^m- \ub^m\right](t_c,x_c)\le\partial_t \ud(t_c,x_c)+\A\ub^m(t_c,x_c)\le \overline{A} \,\kappa_0\,.
\end{equation}
To prove this,
we estimate $\partial_t \ud(t_c,x_c)$ and $\A\ub^m(t_c,x_c)$ separately.
 First we notice that, since $(t_\delta,x_\delta)$ is the first contact point, we  have
\begin{equation}\label{contact.up.1}
\ud(t_\delta,x_\delta)= \ub(t_\delta,x_\delta)\qquad\mbox{and}\qquad
\ud(t,x)\ge \ub(t,x)\,\quad\forall t\in [0,t_\delta]\,,\;\; \forall x\in\Omega\,.
\end{equation}
Hence, since $t_\delta\le t_* $,
\begin{equation}\label{contact.up.2}
\partial_t \ud(t_\delta,x_\delta)\le \partial_t \ub(t_\delta,x_\delta)=\kappa_0\,\Phi_1(x)^{\sigma/m}
\le \kappa_0\,  \|\Phi_1\|_{\LL^\infty(\Omega)}^{\sigma/m} = A_1\,\kappa_0 \,,
\end{equation}
where we defined $A_1:=  \|\Phi_1\|_{\LL^\infty(\Omega)}^{\sigma/m}$. Next we estimate $\A\ub^m(t_c,x_c)$, using the    Kato-type inequality \eqref{Kato.ineq}\,,   namely $\A[u^m]\le mu^{m-1}\A u$\,. This implies
\begin{equation}\label{contact.up.3}\begin{split}
\A[\ub^m](t,x)&\le m\ub^{m-1}(t,x)\, \A \ub(t,x)= m (\kappa_0\, t)^m\,\Phi_1(x)^{\frac{\sigma (m-1)}{m}} \A \Phi_1^\sigma(x)\\
&\le m (\kappa_0\, t_*)^m\,\|\Phi_1\|_{\LL^\infty(\Omega)}^{\frac{\sigma (m-1)}{m}} \left\| \A \Phi_1^\sigma\right\|_{\LL^\infty(\Omega)}
:= A_2\,\kappa_0\,.
\end{split}
\end{equation}
Since $\kappa_0^{m-1}\, t_*^m\leq 1$ (see \eqref{k0.first.cond}),
in order to prove that $A_2$ is finite it is enough to bound $\left\| \A \Phi_1^\sigma\right\|_{\LL^\infty(\Omega)}$.
When $\sigma=1$ we simply have $\A\Phi_1=-\lambda_1\Phi_1$, hence $A_2\leq m \lambda_1 \|\Phi_1\|_{\LL^\infty(\Omega)}^{2-1/m }$. When $\sigma<1$,
we use the assumption $\Phi_1 \in C^\gamma(\Omega)$ to estimate
\begin{equation}\label{contact.up.3a}
|\Phi_1^\sigma(x)-\Phi_1^\sigma(y)|\le  |\Phi_1 (x)-\Phi_1 (y)|^\sigma \le C|x-y|^{\gamma\sigma}\qquad\forall\,x,y\in \Omega\,.
\end{equation}
Hence, since $\gamma\sigma=2sm/(m-1)>2s$ and $K(x,y)\le c_1|x-y|^{-(N+2s)}$, we see that
\[
\begin{split}
\left|\A \Phi_1^\sigma(x)\right|&=\left|\int_{\RR^N}\left[\Phi_1^\sigma(x)-\Phi_1^\sigma(y) \right]K(x,y)\dy\right|\\
&\le \int_{\Omega}|x-y|^{\gamma\sigma} K(x,y)\dy + C\|\Phi_1\|_{\LL^\infty(\Omega)}^\sigma\int_{\RR^N\setminus B_1} |y|^{-(N+2s)}\dy<\infty,
\end{split}
\]
hence $A_2$ is again finite.
Combining \eqref{contact.up.2} and \eqref{contact.up.3}, we obtain \eqref{contact.up} with $\overline{A}:=A_1+A_2$.

\noindent\textit{Lower bound. } We want to prove that there exists $\underline{A}>0$ such that
\begin{equation}\label{contact.low}
-\A\left[u_\delta^m- \ub^m\right](t_c,x_c)\ge\frac{\kb_\Omega}{\|\Phi_1\|_{\LL^\infty(\Omega)}} \int_{\Omega}u_\delta^m(t_c,y)\Phi_1(y) \dy-\underline{A}\,\kappa_0\,.
\end{equation}
This follows by (L1)
and  \eqref{contact.up.1} as follows:
\begin{equation}\label{contact.4aaaaa}\begin{split}
-\A &\left[u_\delta^m-\ub^m\right](t_c,x_c)\\
&=-\int_{\RR^N}\left[\big(u_\delta^m(t_c,x_c)-u_\delta^m(t_c,y)\big)-\big(\ub^m(t_c,x_c)-\ub^m(t_c,y)\big)\right]K(x,y)\dy\\
&=\int_{\Omega}\left[u_\delta^m(t_c,y)-\ub^m(t_c,y)\right]K(x,y)\dy\\
&\ge \kb_\Omega \int_{\Omega}\left[u_\delta^m(t_c,y)-\ub^m(t_c,y)\right] \dy
\ge \frac{\kb_\Omega}{\|\Phi_1\|_{\LL^\infty(\Omega)}} \int_{\Omega}u_\delta^m(t_c,y)\Phi_1(y) \dy-\underline{A}\,\kappa_0,
\end{split}\end{equation}
where in the last step we used that $\ub^m(t_c,y)= [\k_0 t \Phi_1^{\sigma/m}(y)]^m   \le \k_2(\k_0 t_*)^m $ and $\k_0^{m-1}t_*^m\le 1$ (see \eqref{k0.first.cond}).

\noindent\textit{End of the proof.} The contradiction can be now obtained by joining the upper and lower bounds \eqref{contact.up} and \eqref{contact.low}.
More precisely, we have proved that
\[
\int_{\Omega}u_\delta^m(t_c,y)\Phi_1(y) \dy\le \frac{\|\Phi_1\|_{\LL^\infty(\Omega)}}{\kb_\Omega}(\overline{A}+\underline{A})\,\kappa_0 :=\ka\,\k_0,
\]
that combined with the lower bound \eqref{lem1.lower.bdd} yields
$$
c_2 \left(\int_{\Omega}u_0(x)\Phi_1(x)\dx \right)^m\le  \int_{\Omega}u_\delta^m(t_c,y)\Phi_1(y)\dy
\le \ka\,\k_0.
$$
Setting $\k_0 := \left(1\wedge \frac{c_2}{\ka}\right)t_*^{-m/(m-1)}$ we obtain the desired contradiction.\qed

\subsection{Proof of Theorem \ref{thm.Lower.PME.large.t}.~}
We first recall the upper pointwise estimates \eqref{thm.NLE.PME.estim}: for all $0\le t_0\le t_1 \le t $ and a.e. $x_0\in \Omega$\,, we have that
\begin{equation}\label{Lower.PME.Step.1.1}
\int_{\Omega}u(t_0,x)\K(x , x_0)\dx - \int_{\Omega}u({t_1},x)\K(x , x_0)\dx \le (m-1)\frac{t^{\frac{m}{m-1}}}{t_0^{\frac{1}{m-1}}} \,u^m(t,x_0)\,.
\end{equation}
The proof follows by estimating the two integrals on the left-hand side separately.

We begin by using the upper bounds \eqref{thm.Upper.PME.Boundary.2} to get
\begin{equation}\label{Lower.PME.Step.1.2}
\int_{\Omega}u(t_1,x)\K(x, x_0)\dx\le \ka \frac{\Phi_1(x_0)}{t_1^{\frac{1}{m-1}}}\qquad\mbox{for all $(t_1,x)\in(0,+\infty)\times\Omega$}\,.
\end{equation}
Then we note that, as a consequence of (K2) and Lemma \ref{cor.abs.L1.phi},
\begin{equation}\label{Lower.PME.Step.1.4}
\int_{\Omega}u(t_0,x)\K(x, x_0)\dx\ge \kb_\Omega \Phi_1(x_0)\int_{\Omega}u(t_0,x)\Phi_1(x)\dx
\ge \frac{\kb_\Omega}{2}\Phi_1(x_0)\int_{\Omega}u_0(x)\Phi_1(x)\dx
\end{equation}
provided $t_0\le \frac{\tau_0}{\|u_0\|_{\LL^1_{\Phi_1}(\Omega)}^{m-1}}$.
Combining \eqref{Lower.PME.Step.1.1}, \eqref{Lower.PME.Step.1.2}, and \eqref{Lower.PME.Step.1.4}, for all $t\ge t_1\ge t_0\ge 0$ we obtain
$$
u^m(t,x_0)
\ge \frac{t_0^{\frac{1}{m-1}}}{m-1}\left(\frac{\kb_\Omega}{2}\|u_0\|_{\LL^1_{\Phi_1}(\Omega)}
 -  \ka t_1^{-\frac{1}{m-1}}\right)\frac{\Phi_1(x_0)}{t^{\frac{m}{m-1}}}\,.
$$
Choosing
\[
t_0:=\frac{\tau_0}{\|u_0\|_{\LL^1_{\Phi_1}(\Omega)}^{m-1}}\le t_1:=t_*=\frac{\k_*}{\|u_0\|_{\LL^1_{\Phi_1}(\Omega)}^{m-1}}
\qquad\mbox{with}\qquad
\k_*\ge \tau_0 \vee \left(\frac{\kb_\Omega}{4\ka}\right)^{m-1}
\]
so that
$
\frac{\kb_\Omega}{2}\|u_0\|_{\LL^1_{\Phi_1}(\Omega)} -  \ka t_1^{-\frac{1}{m-1}}
\ge \frac{\kb_\Omega}{4} \|u_0\|_{\LL^1_{\Phi_1}(\Omega)},$
the result follows.\qed

\subsection{Proof of Theorems \ref{prop.counterex} and \ref{prop.counterex2}}

\noindent\textbf{Proof of Theorem \ref{prop.counterex}. }
Since $u_0\leq C_0 \Phi_1$ and $\A \Phi=\lambda_1 \Phi_1$, we have
$$
\int_\Omega u_0(x) \K(x,x_0)\dx\leq C_0\int_\Omega \Phi_1(x)\K(x,x_0)\dx=C_0 \AI\Phi_1(x_0)=\frac{C_0}{\lambda_1}\Phi_1(x_0).
$$
Since $t\mapsto \int_\Omega u(t,y)\K(x,y)\dy$ is decreasing (see \eqref{thm.NLE.PME.estim.0}), it follows  that
\begin{equation}
\label{eq:initial data.2}
\int_\Omega u(t,y)\K(x_0,y)\dy\le \frac{C_0}{\lambda_1}\Phi_1(x_0)\qquad\mbox{for all }t\ge 0.
\end{equation}
Combining this estimate with \eqref{Upper.PME.Step.1.2.c} concludes the proof.
\qed

\noindent\textbf{Proof of Theorem \ref{prop.counterex2}. }
Given $x_0 \in \Omega$, set $R_0:={\rm dist}(x_0,\partial\Omega)$.
Since $\K(x,x_0) \gtrsim |x-x_0|^{-(N-2s)}$ inside $B_{R_0/2}(x_0)$ (by (K4)),
using our assumption on $u(T)$ we get
$$
\int_\Omega \K(x,x_0) u(T,x)\dx\gtrsim
\int_{B_{R_0/2}(x_0)} \frac{  \Phi_1(x)^\alpha}{|x-x_0|^{N-2s}}
\gtrsim \Phi_1(x_0)^\alpha R_0^{2s}.
$$
Recalling that $\Phi_1(x_0)\asymp R_0^\gamma$, this yields
$$
\Phi_1(x_0)^{\alpha+\frac{2s}{\gamma}}\lesssim \int_\Omega \K(x,x_0) u(T,x)\dx.
$$
Combining the above inequality with \eqref{eq:initial data.2} gives
\[
\Phi_1(x_0)^{\alpha+\frac{2s}{\gamma}} \lesssim  \Phi_1(x_0) \qquad \forall\,x_0\in \Omega,
\qquad\mbox{which implies}\qquad \alpha \geq 1-\frac{2s}{\gamma}\,.
\]
Noticing that $1-\frac{2s}{\gamma}>\frac{1}{m}$ if and only if $\sigma<1$, this concludes the proof.\qed

\section{Summary of the general decay and boundary results}\label{sect.Harnack}

In this section we present as summary of the main results, which can be summarized in various forms of upper and lower bounds, which we call Global Harnack Principle, (GHP) for short.  As already mentioned, such inequalities are important for regularity issues (see Section \ref{sect.regularity}), and  they play a fundamental role in formulating the sharp asymptotic behavior (see Section \ref{sec.asymptotic}). The proof of such GHP is obtained by combining upper and lower bounds, stated and proved in Sections \ref{sect.upperestimates} and \ref{Sec.Lower} respectively. There are  cases when the bounds do not match, for which the complicated panorama  described in the Introduction holds. As explained before, as far as examples are concerned, the latter anomalous situation happens only for the SFL.

\begin{thm}[Global Harnack Principle I]\label{thm.GHP.PME.I}
Let $\A$ satisfy (A1), (A2), (K2), and (L1). Furthermore, suppose that $\A$ has a first eigenfunction $\Phi_1\asymp \dist(x,\partial\Omega)^\gamma$\,. Let $\sigma$ be as in \eqref{as.sep.var} and assume that $2sm\ne \gamma(m-1)$ and:\\
- either $\sigma=1$;\\
- or $\sigma<1$, $K(x,y)\le c_1|x-y|^{-(N+2s)}$ for a.e. $x,y\in \RR^N$,
and $\Phi_1\in C^\gamma(\overline{\Omega})$.\\
Let $u\ge 0$ be a weak dual solution to the (CDP) corresponding to $u_0\in \LL^1_{\Phi_1}(\Omega)$. Then, there exist  constants $\kb,\ka>0$, so that the following inequality holds:
\begin{equation}\label{thm.GHP.PME.I.Ineq}
\kb\, \left(1\wedge \frac{t}{t_*}\right)^{\frac{m}{m-1}}\frac{\Phi_1(x)^{\sigma/m}}{t^{\frac{1}{m-1}}}
\le \, u(t,x) \le \ka\, \frac{\Phi_1(x)^{\sigma/m}}{t^{\frac1{m-1}}}\qquad\mbox{for all $t>0$ and all $x\in \Omega$}\,.
\end{equation}
The constants $\kb,\ka$  depend only on   $N,s,\gamma, m, c_1,\kb_\Omega,\Omega$, and $\|\Phi_1\|_{C^\gamma(\Omega)}$\,.
\end{thm}
\noindent {\sl Proof.~}We combine the upper bound \eqref{thm.Upper.PME.Boundary.F}
with the lower bound \eqref{Thm.lower.B}. The expression of $t_*$ is explicitly given in Theorem \ref{Thm.lower.B}. \qed
\medskip
\noindent\textbf{Degenerate   kernels. }When the kernel $K$ vanishes on $\partial\Omega$, there are two  combinations of upper/lower bounds that provide Harnack inequalities, one for small times and one for large times.
As we have already seen, there is a strong difference between the case $\sigma=1$ and $\sigma<1$.

\begin{thm}[Global Harnack Principle II]\label{thm.GHP.PME.II}
Let  (A1), (A2), and (K2) hold. Let $u\ge 0$ be a weak dual solution to the (CDP) corresponding to $u_0\in \LL^1_{\Phi_1}(\Omega)$. Assume that:\\
- either $\sigma=1$ and $2sm\ne \gamma(m-1)$;\\
- or $\sigma<1$, $u_0\geq \kb_0\Phi_1^{\sigma/m}$ for some $\kb_0>0$, and (K4) holds.\\
Then there exist  constants $\kb,\ka>0$ such that the following inequality holds:
$$
\kb\,  \frac{\Phi_1(x)^{\sigma/m}}{t^{\frac{1}{m-1}}}
\le \, u(t,x) \le \ka\, \frac{\Phi_1(x)^{\sigma/m}}{t^{\frac1{m-1}}}\qquad\mbox{for all $t\ge t_*$ and all $x\in \Omega$}\,.
$$
When $2sm= \gamma(m-1)$,  assuming (K4)
and that $u_0\geq \kb_0\Phi_1\left(1+|\log\Phi_1 |\right)^{1/(m-1)}$ for some $\kb_0>0$, then for all $t\ge t_*$ and all $x\in \Omega$
\begin{multline*}
\kb\,  \frac{\Phi_1(x)^{1/m}}{t^{\frac{1}{m-1}}}\left(1+|\log\Phi_1(x) |\right)^{1/(m-1)}
\le \, u(t,x)
\le \ka\, \frac{\Phi_1(x)^{1/m} }{t^{\frac1{m-1}}}\left(1+|\log\Phi_1(x) |\right)^{1/(m-1)}\,.
\end{multline*}
The constants $\kb,\ka$ depend only on   $N,s,\gamma, m, \kb_0,\kb_\Omega$, and $\Omega$.
\end{thm}
\noindent {\sl Proof.~}In the case $\sigma=1$, we combine the upper bound \eqref{thm.Upper.PME.Boundary.F} with the lower bound \eqref{thm.Lower.PME.Boundary.large.t}. The expression of $t_*$ is explicitly given in Theorem \ref{thm.Lower.PME.large.t}.
When $\sigma<1$, the upper bound is still given \eqref{thm.Upper.PME.Boundary.F},
while the lower bound follows by comparison with the solution $S(x)(\kb_0^{1-m}+t)^{-\frac{1}{m-1}}$, recalling that $S\asymp \Phi_1^{\sigma/m}$ (see Theorem \ref{Thm.Elliptic.Harnack.m}).
\qed

\noindent\textbf{Remark. }Local Harnack inequalities of elliptic/backward type follow as a consequence of Theorems \ref{thm.GHP.PME.I} and \ref{thm.GHP.PME.II}, for all times and for large times respectively, see Theorem \ref{thm.Harnack.Local}.

\medskip

Note that, for small times, we cannot find matching powers for a global Harnack inequality (except for some special initial data), and such result is actually false for $s=1$ (in view of the finite speed of propagation).
Hence, in the remaining cases, we have only the following general result.
\begin{thm}[{Non matching upper and lower bounds}]\label{thm.GHP.PME.III}
Let $\A$ satisfy (A1), (A2), (K2), and (L2).
Let $u\ge 0$ be a weak dual solution to the (CDP) corresponding to $u_0\in \LL^1_{\Phi_1}(\Omega)$. Then, there exist  constants $\kb,\ka>0$, so that the following inequality holds  when $2sm\ne \gamma(m-1)$:
\begin{equation}\label{thm.GHP.PME.III.Ineq}
\kb\,\left(1\wedge \frac{t}{t_*}\right)^{\frac{m}{m-1}}\frac{\Phi_1(x)}{t^{\frac{1}{m-1}}}
\le \, u(t,x) \le \ka\, \frac{\Phi_1(x)^{\sigma/m}}{t^{\frac1{m-1}}}\qquad\mbox{for all $t>0$ and all $x\in \Omega$}.
\end{equation}
When $2sm=\gamma(m-1)$, a logarithmic correction $\left(1+|\log\Phi_1(x) |\right)^{1/(m-1)}$ appears in the right hand side.
\end{thm}
\noindent {\sl Proof.~} We combine the upper bound \eqref{thm.Upper.PME.Boundary.F} with the lower bound \eqref{thm.Lower.PME.Boundary.1}. The expression of $t_*$ is explicitly given in Theorem \ref{thm.Lower.PME}.\qed

\noindent\textbf{Remark. }As already mentioned in the introduction, in the non-matching case, which in examples can only happen for spectral-type operators, we have the appearance of an \textit{anomalous behaviour of solutions }corresponding to ``small data'': it happens for all times when $\sigma<1$ or $2sm=\gamma(m-1)$, and it can eventually happen for short times when $\sigma=1$.
%

\section{Asymptotic behavior}\label{sec.asymptotic}

An important application of the Global Harnack inequalities of the previous section concerns the sharp asymptotic behavior of solutions. More precisely, we first show that for large times all solutions behave like the separate-variables solution $\mathcal{U}(t,x)={S(x)}\,{t^{-\frac{1}{m-1}}}$
introduced at the end of Section \ref{sec.results}.
Then, whenever the (GHP) holds, we can improve this result to an estimate in relative error.

\begin{thm}[Asymptotic behavior]\label{Thm.Asympt.0}
Assume that $\A$ satisfies (A1), (A2), and (K2), and let $S$ be as in Theorem \ref{Thm.Elliptic.Harnack.m}. Let $u$ be any weak dual solution to the (CDP).
Then, unless $u\equiv 0$,
\begin{equation}\label{Thm.Asympt.0.1}
\left\|t^{\frac{1}{m-1}}u(t,\cdot)- S\right\|_{\LL^\infty(\Omega)}\xrightarrow{t\to\infty}0\,.
\end{equation}
\end{thm}
\noindent {\sl Proof.~}The proof uses rescaling and time monotonicity arguments, and it is a simple adaptation of the proof of Theorem 2.3 of  \cite{BSV2013}. In those arguments, the interior $C^{\alpha}_{x}(\Omega)$ continuity is needed to improve the $\LL^1(\Omega)$ convergence to $\LL^\infty(\Omega)$, but the interior H\"older continuity is guaranteed
by Theorem \ref{thm.regularity.1}(i) below. \qed

We now exploit the (GHP) to get a stronger result.

\begin{thm}[Sharp asymptotic behavior]\label{Thm.Asympt}
Under the assumptions of Theorem \ref{Thm.Asympt.0}, assume that $u\not\equiv 0$.
Furthermore,
suppose that either the assumptions of Theorem \ref{thm.GHP.PME.I} or of Theorem \ref{thm.GHP.PME.II} hold.
Set $\mathcal{U}(t,x):=t^{-\frac{1}{m-1}}S(x)$.
Then there exists $c_0>0$ such that, for all $t\ge t_0:=c_0\|u_0\|_{\LL^1_{\Phi_1}(\Omega)}^{-(m-1)}$, we have
\begin{equation}\label{conv.rates.rel.err}
\left\|\frac{u(t,\cdot)}{\mathcal{U}(t,\cdot)}-1 \right\|_{\LL^\infty(\Omega)} \le \frac{2}{m-1}\,\frac{t_0}{t_0+t}\,.
\end{equation}
We remark that the constant $c_0>0$ only depends on $N,s,\gamma, m, \kb_0,\kb_\Omega$, and $\Omega$.
\end{thm}
\noindent\textbf{Remark. }This asymptotic result is sharp, as it can be checked by considering $u(t,x)=\mathcal{U}(t+1,x)$.
For the classical case, that is $\A=\Delta$, we recover the classical results of \cite{Ar-Pe,JLVmonats} with a different proof.

\noindent {\sl Proof.~}Notice that we are in the position to use of Theorems \ref{thm.GHP.PME.I} or \ref{thm.GHP.PME.II}, namely we have
$$
u(t) \asymp t^{-\frac1{m-1}}S=\mathcal{U}(t,\cdot) \qquad\mbox{for all $t\ge t_*$}\,,
$$
where the last equivalence follows by Theorem \ref{Thm.Elliptic.Harnack.m}. Hence, we can rewrite the bounds above saying that there exist $\kb,\ka>0$ such that
\begin{equation}\label{GHP.ASYM}
\kb\, \frac{S(x)}{t^{\frac1{m-1}}}
\le \, u(t,x) \le \ka\, \frac{S(x)}{t^{\frac1{m-1}}}\qquad\mbox{for all $t\ge t_*$ and a.e. $x\in \Omega$}\,.
\end{equation}
Since $t_*=\kappa_*\|u_0\|_{\LL^1_{\Phi_1}(\Omega)}^{-(m-1)}$,
the first inequality implies that
$$
\frac{S}{(t_*+t_0)^{\frac{1}{m-1}}}\leq \kb\, \frac{S}{t_*^{\frac1{m-1}}} \leq u(t_*)
$$
for some $t_0=c_0\|u_0\|_{\LL^1_{\Phi_1}(\Omega)}^{-(m-1)} \geq t_*$.
Hence, by comparison principle,
$$
\frac{S}{(t+t_0)^{\frac{1}{m-1}}}
\le \, u(t)\qquad\mbox{for all $t\ge t_*$.}
$$
On the other hand, it follows by \eqref{GHP.ASYM} that $u(t,x) \leq \mathcal U_T(t,x):=S(x)(t-T)^{-\frac1{m-1}}$ for all $t \geq T$ provided $T$ is large enough.
If we now start to reduce $T$, the comparison principle combined with the upper bound
\eqref{thm.Upper.PME.Boundary.F} shows that $u$ can never touch $\mathcal U_T$ from below in $(T,\infty)\times \Omega$. Hence we can reduce $T$ until $T=0$, proving that $u \leq \mathcal U_0$ (for an alternative proof, see Lemma 5.4 in
\cite{BSV2013}).
Since $t_0 \geq t_*$,
this shows that
$$
\frac{S(x)}{(t+t_0)^{\frac{1}{m-1}}}
\le \, u(t,x) \le \frac{S(x)}{t^{\frac1{m-1}}}\qquad\mbox{for all $t\ge t_0$ and a.e. $x\in \Omega$}\,,
$$
therefore
$$
\biggl|1-\frac{u(t,x)}{\mathcal{U}(t,x)}\biggr| \leq 1-\biggl(1-\frac{t_0}{t_0+t}\biggr)^{\frac{1}{m-1}}\leq \frac{2}{m-1}\frac{t_0}{t_0+t}\qquad\mbox{for all $t\ge t_0$ and a.e. $x\in \Omega$},
$$
as desired.\qed


\section{Regularity results}\label{sect.regularity}

\noindent In order to obtain the regularity results, we basically require the validity of a Global Harnack Principle, namely Theorems \ref{thm.GHP.PME.I}, \ref{thm.GHP.PME.II}, or \ref{thm.GHP.PME.III}, depending on the situation under study. For some higher regularity results, we will eventually need some extra assumptions on the kernels.
For simplicity we assume that $\A$ is described by a kernel, without any lower order term. However, it is clear that the presence of lower order terms does not play any role in the interior regularity. 

\begin{thm}[Interior Regularity]\label{thm.regularity.1}
Assume that
\begin{equation*}
\A f(x)=P.V.\int_{\RR^N}\big(f(x)-f(y)\big)K(x,y)\dy+B(x)f(x)\,,
\end{equation*}
with
$$
K(x,y)\asymp |x-y|^{-(N+2s)}\quad \text{ in $B_{2r}(x_0)\subset\Omega$},
\qquad
K(x,y)\lesssim |x-y|^{-(N+2s)}\quad \text{ in $\mathbb R^N\setminus B_{2r}(x_0)$}.
$$
Let $u$ be a nonnegative bounded weak dual solution to problem (CDP) on $(T_0,T_1)\times\Omega$, and assume that there exist $\delta,M>0$ such that
\begin{align*}
&0<\delta\le u(t,x)\qquad\mbox{for a.e. $(t,x)\in (T_0,T_1)\times B_{2r}(x_0)$,}\\
&0\leq u(t,x)\leq M\qquad\mbox{for a.e. $(t,x)\in (T_0,T_1)\times \Omega$.}
\end{align*}\vspace{-7mm}
\begin{enumerate}
\item[(i)] Then $u$ is {H\"older continuous in the interior}. More precisely, there exists $\alpha>0$ such that, for all $0<T_0<T_2<T_1$,
\begin{equation}\label{thm.regularity.1.bounds.0}
\|u\|_{C^{\alpha/2s,\alpha}_{t,x}((T_2,T_1)\times B_{r}(x_0))}\leq C.
\end{equation}
\item[(ii)] Assume in addition $|K(x,y)-K(x',y)|\le c |x-x'|^\beta\,|y|^{-(N+2s)}$ for some $\beta\in (0,1\wedge 2s)$ such that $\beta+2s$ is not an integer. Then {$u$ is a classical solution in the interior}. More precisely, for all $0<T_0<T_2<T_1$,
\begin{equation}\label{thm.regularity.1.bounds.1}
\|u\|_{C^{1+\beta/2s, 2s+\beta}_{t,x}((T_2,T_1)\times  B_{r}(x_0))}\le C.
\end{equation}
\end{enumerate}
\end{thm}
The constants in the above regularity estimates depend on the solution only through the upper and lower bounds on $u$. This bounds can be made quantitative by means of local Harnack inequalities, of elliptic and forward type, which follows from the Global ones.
\begin{thm}[Local Harnack Inequalities of Elliptic/Backward Type]\label{thm.Harnack.Local}
Under the assumptions of Theorem \ref{thm.GHP.PME.I},
there exists a constant $\hat H>0$, depending only on $N,s,\gamma, m, c_1,\kb_\Omega,\Omega$, such that for all balls $B_R(x_0)$ such that $B_{2R}(x_0)\subset\Omega$:
\begin{equation}\label{thm.Harnack.Local.ell}
\sup_{x\in B_R(x_0)}u(t,x)\le \frac{\hat H}{\left(1\wedge\frac{t}{t_*}\right)^{\frac{m}{m-1}}}\, \inf_{x\in B_R(x_0)}u(t,x)\qquad\mbox{for all }t>0.
\end{equation}
Moreover, for all $t>0$ and all $h>0$ we have the following:
\begin{equation}\label{thm.Harnack.Local.ForwBack}
\sup_{x\in B_R(x_0)}u(t,x)\le \hat H\left[\left(1+\frac{h}{t}\right) \left(1\wedge\frac{t}{t_*}\right)^{-m}\right]^{\frac{1}{m-1}}\, \inf_{x\in B_R(x_0)}u(t+h,x)\,.
\end{equation}
\end{thm}
\noindent {\sl Proof.~}
Recalling \eqref{thm.GHP.PME.I.Ineq}, the bound \eqref{thm.Harnack.Local.ell}
follows easily from the following Harnack inequality for the first eigenfunction, see for instance \cite{BFV-Elliptic}:
$$
\sup_{x\in B_R(x_0)} \Phi_1(x)\le H_{N,s,\gamma,\Omega} \inf_{x\in B_R(x_0)} \Phi_1(x).
$$
Since $u(t,x)\leq (1+h/t)^{\frac{1}{m-1}}u(t+h,x)$, by the time monotonicity of $t\mapsto t^{\frac{1}{m-1}}\,u(t,x)$,
\eqref{thm.Harnack.Local.ForwBack} follows.\qed

\noindent\textbf{Remark. }The same result holds for large times, $t\ge t_*$ as a consequence of Theorem \ref{thm.GHP.PME.II}. Already in the local case $s=1$\,, these Harnack inequalities are stronger than the known two-sided inequalities valid for solutions to the Dirichlet problem for the classical porous medium equation,  cf. \cite{AC83,DaskaBook,DiB88,DiBook,DGVbook}, which are of forward type and are often stated in terms of the so-called intrinsic geometry. Note that elliptic and backward Harnack-type inequalities usually occur in the fast diffusion range $m<1$ \cite{BGV-Domains, BV, BV-ADV, BV2012}, or for linear equations in bounded domains  \cite{FGS86,SY}.

For sharp boundary regularity we need a GHP with matching powers, like Theorems \ref{thm.GHP.PME.I} or \ref{thm.GHP.PME.II}, and when $s>\gamma/2$, we can also prove H\"older regularity up to the boundary. We leave to the interested reader to check that the presence of an extra term $B(x)u^m(t,x)$ with $0\leq B(x)\leq c_1 {\rm dist}(x,\partial\Omega)^{-2s}$ (as in the SFL) does not affect the validity of the next result. Indeed, when considering the scaling in \eqref{eq:scaling}, the lower term scales as $\hat B_r u_r^m$ with $0\leq \hat B_r\leq c_1$ inside the unit ball $B_1$.
\begin{thm}[H\"older continuity up to the boundary]\label{thm.regularity.2}
Under  assumptions of Theorem $\ref{thm.regularity.1}$(ii),
assume in addition that $2s>\gamma$
Then {$u$ is H\"older continuous up to the boundary}. More precisely, for all $0<T_0<T_2<T_1$ there exists  a constant $C>0$ such that
\begin{equation}\label{thm.regularity.2.bounds}
\|u\|_{C^{\frac{\gamma}{m\vartheta},\frac{\gamma}{m}}_{t,x}((T_2,T_1)\times{\Omega})}\le C
\qquad\mbox{with }\quad\vartheta:=2s-\gamma\left(1-\frac{1}{m}\right)\,.
\end{equation}
\end{thm}

\noindent\textbf{Remark. } Since we have $u(t,x)\asymp \Phi_1(x)^{1/m}\asymp\dist(x,\partial\Omega)^{\gamma/m}$ (note that $2s>\gamma$ implies that $\sigma=1$  and that $2sm\ne \gamma(m-1)$),  the spacial H\"older exponent is sharp, while the H\"older exponent in time is the natural one by scaling.

\subsection{Proof of interior regularity}

The strategy to prove Theorem \ref{thm.regularity.1} follows the lines of \cite{BFR} but with some modifications.
The basic idea is that, because $u$ is bounded away from zero and infinity, the equation is non-degenerate and we can use parabolic regularity for nonlocal equations to obtain the results. More precisely, interior H\"older regularity will follow by applying $C^{\alpha/2s,\alpha}_{t,x}$ estimates of \cite{FK} for a ``localized'' linear problem.
Once H\"older regularity is established, under an H\"older continuity assumption on the kernel we can use the Schauder estimates proved in \cite{Schauder} to conclude.

\subsubsection{Localization of the problem }\label{sssect.localization}
Up to a rescaling, we can assume $r=2$, $T_0=0$, $T_1=1$. Also, by a standard covering argument,
it is enough to prove the results with $T_2=1/2$.

Take a cutoff function $\rho\in C^\infty_c(B_4)$ such that $\rho\equiv 1$ on $B_3$,
$\eta\in C^\infty_c(B_2)$ a cutoff function such that $\eta\equiv 1$ on $B_1$,
 and define $v=\rho u$. By construction $u=v$ on $(0,1)\times B_3$. Since $\rho\equiv 1$ on $B_3$, we can write the equation for $v$ on the small cylinder $(0,1)\times B_1$ as
$$
\partial_t v(t,x) = -\A [v^m](t,x) +g(t,x)= - L_a v(t,x) + f(t,x) +g(t,x)
$$
where
$$
L_a[v](t,x):=\int_{\RR^N}\big(v(t,x)-v(t,y)\big)a(t,x,y)K(x,y)\dy\,,
$$
\begin{align*}
a(t,x,y) &:= \frac{v^m(t,x)-v^m(t,y)}{v(t,x)-v(t,y)}\eta(x-y)+ \big[1-\eta(x-y)\big]\\
&=m\eta(x-y)\int_0^1 \left[(1-\lambda)v(t,x) +\lambda v(t,y) \right]^{m-1} \rd\lambda   + \big[1-\eta(x-y)\big]\,,
\end{align*}
\[
f(t,x):=\int_{\RR^N\setminus B_1(x)}\bigl(v^m(t,x)-v^m(t,y)-v(t,x)+v(t,y)\bigr)[1-\eta(x-y)] K(x,y)\dy\,,
\]
and
$$
g(t,x):= -\A\left[(1-\rho^m)u^m\right](t,x)=\int_{\RR^N\setminus B_3}(1-\rho^m(y))u^m(t,y)  K(x,y)\dy
$$
(recall that $(1-\rho^m)u^m\equiv 0$ on $(0,1)\times B_3$).

\subsubsection{H\"older continuity in the interior}\label{sssect.Calpha}

Set $b:=f+g$, with $f$ and $g$ as above. It is easy to check that, since $K(x,y)\lesssim |x-y|^{-(N+2s)}$, $b \in \LL^\infty((0,1)\times B_1)$.
Also, since $0<\delta\leq u \leq M$ inside $(0,1)\times B_1$, there exists $\Lambda>1$ such that
$\Lambda^{-1}\leq a(t,x,y)\leq \Lambda$ for a.e. $(t,x,y)\in (0,1)\times B_1\times B_1$ with $|x-y|\leq 1$.
This guarantees that the linear operator $L_a$ is uniformly elliptic,
so we can apply the results in \cite{FK} to ensure that
\[\|v\|_{C^{\alpha/2s,\alpha}_{t,x}((1/2,1)\times B_{1/2})}\leq C\bigl(\|b\|_{L^\infty((0,1)\times B_1)}+\|v\|_{L^\infty((0,1)\times\RR^N)}\bigr)\]
for some universal exponent $\alpha>0$. This proves Theorem \ref{thm.regularity.1}(i).

\subsubsection{Classical solutions in the interior}
Now that we know that $u\in C^{\alpha/2s,\alpha}((1/2,1)\times B_{1/2})$,
we repeat the localization argument above with cutoff functions $\rho$ and $\eta$ supported inside $(1/2,1)\times B_{1/2}$ to ensure that $v:=\rho u$ is H\"older continuous in $(1/2,1)\times\RR^N$. Then,
to obtain higher regularity we argue as follows.

Set $\beta_1:=\min\{\alpha,\beta\}$.
Thanks to the assumption on $K$ and Theorem \ref{thm.regularity.1}(i), it is easy to check that $K_a(t,x,y):=a(t,x,y)K(x,y)$
satisfies
$$
|K_a(t,x,y)-K_a(t',x',y)| \leq C\left(|x-x'|^{\beta_1}+|t-t'|^{\beta_1/2s}\right)|y|^{-(N+2s)}
$$
inside $(1/2,1)\times B_{1/2}$.
Also, $f,g \in C^{\beta_1/2s,\beta_1}((1/2,1)\times B_{1/2})$.
This allows us to apply the Schauder estimates from \cite{Schauder} (see also \cite{CKriv})
to obtain that
$$
\|v\|_{C^{1+\beta_1/2s,2s+\beta}_{t,x}((3/4,1)\times B_{1/4})}
\leq C\bigl(\|b\|_{C^{\beta/2s,\beta}_{t,x}((1/2,1)\times B_{1/2})}+\|v\|_{{C^{\beta/2s,\beta}_{t,x}((1/2,1)\times\RR^N)}}\bigr).
$$
In particular, $u\in C^{1+\beta_1/2s,2s+\beta_1}((3/4,1)\times B_{1/8})$.
In case $\beta_1=\beta$ we stop here. Otherwise we set
$\alpha_1:=2s+\beta$ and we repeat the argument above with $\beta_2:=\min\{\alpha_1,\beta\}$
in place of $\beta_1$. In this way, we obtain that $u\in C^{1+\beta_1/2s,2s+\beta_1}((1-2^{-4},1)\times B_{2^{-5}})$. Iterating this procedure finitely many times, we finally obtain that
$$
u\in C^{1+\beta/2s,2s+\beta}((1-2^{-k},1)\times B_{2^{-k-1}})
$$
for some universal $k$. Finally, a covering argument completes the proof
of Theorem \ref{thm.regularity.1}(ii).
\subsection{Proof of boundary regularity}
The proof of Theorem \ref{thm.regularity.2} follows by scaling and interior estimates.
Notice that the assumption $2s>\gamma$ implies that $\sigma=1$,
hence $u(t)$ has matching upper and lower bounds.

Given $x_0 \in \Omega$, set $r=\textrm{dist}(x_0,\partial\Omega)/2$
and define
\begin{equation}
\label{eq:scaling}
u_r(t,x):=r^{-\frac{\gamma}{m}}\,u\left(t_0+r^\vartheta t,\,x_0+rx\right)\,, \qquad\mbox{with}\qquad\vartheta:=2s-\gamma\left(1-\frac{1}{m}\right)\,.
\end{equation}
Note that, because $2s>\gamma$, we have $\vartheta>0$.
With this definition, we see that $u_r$ satisfies the equation $\partial_t u_r+\A_r u_r^m=0$ in $\Omega_r:=(\Omega-x_0)/r$, where
\[
\A_r f(x)=P.V.\int_{\RR^N}\big(f(x)-f(y)\big)K_r(x,y)\dy\,,\qquad K_r(x,y):=r^{N+2s}K(x_0+rx,x_0+ry)\,.
\]
Note that, since $\sigma=1$, it follows by the (GHP) that $u(t)\asymp \dist(x,\partial\Omega)^{\gamma/m}$. Hence,
\[
0<\delta\leq u_r(t,x)\leq M, \qquad\mbox{for all $t\in[r^{-\vartheta}T_0,r^{-\vartheta}T_1]$ and $x\in B_1$,}
\]
with constants $\delta,M>0$ that are independent of $r$ and $x_0$.
In addition, using again that $u(t)\asymp \dist(x,\partial\Omega)^{\gamma/m}$, we see that
\[
u_r(t,x)\leq C\bigl(1+|x|^{\gamma/m}\bigr)
\qquad \text{for all $t\in[r^{-\vartheta}T_0,r^{-\vartheta}T_1]$ and $x\in \RR^N$}.
\]
Noticing that
$
u_r^m(t,x) \leq C(1+|x|^{\gamma})
$
and that $\gamma<2s$ by assumption, we see that the tails
of $u_r$ will not create any problem. Indeed, for any $x \in B_1$,
$$
\int_{\mathbb R^N\setminus B_2}u_r^m(t,y)K_r(x,y)^{-(N+2s)}\,dy\leq C\int_{\mathbb R^N\setminus B_2}|y|^\gamma |y|^{-(N+2s)}\,dy\leq \bar C_0,
$$
where $\bar C_0$ is independent of $r$.
This means that we can localize the problem using cutoff functions as done in Section \ref{sssect.localization}, and the integrals defining the functions $f$ and $g$
will converge uniformly with respect to $x_0$ and $r$.
Hence, we can apply Theorem \ref{thm.regularity.1}(ii) to get
\begin{equation}
\label{eq:interior ur}
\|u_r\|_{C^{1+\beta/2s,2s+\beta}([r^{-\vartheta}T+1/2,r^{-\vartheta}T+1]\times B_{1/2})} \leq C
 \end{equation}
for all $T \in [T_0,T_1-r^{-\vartheta}]$.
Since $\gamma/m<2s+\beta$ (because $\gamma<2s$),
it follows that
$$
\|u_r\|_{L^\infty([r^{-\vartheta}T+1/2,r^{-\vartheta}T+1],C^{\gamma/m}(B_{1/2})} \leq \|u_r\|_{C^{1+\beta/2s,2s+\beta}([r^{-\vartheta}T+1/2,r^{-\vartheta}T_0+1]\times B_{1/2})} \leq C.
$$
Noticing that
$$
\sup_{t \in [r^{-\vartheta}T+1/2,r^{-\vartheta}T+1]}[u_r]_{C^{\gamma/m}(B_{1/2})}
=
\sup_{t \in [T+r^{\vartheta}/2,r^{-\vartheta}T+r^{-\vartheta}]}[u]_{C^{\gamma/m}(B_{r}(x_0))},
$$
and that $T \in [T_0,T_1-r^{-\vartheta}]$ and $x_0$ are arbitrary, arguing as in \cite{RosSer}
we deduce that, given $T_2\in (T_0,T_1)$,
\begin{equation}
\label{eq:holder x}
\sup_{t \in [T_2,T_1]}[u]_{C^{\gamma/m}(\Omega)}
 \leq C.
 \end{equation}
This proves the global H\"older regularity  in space.
To show the regularity in time, we start again from \eqref{eq:interior ur} to get
$$
\|\partial_tu_r\|_{\LL^\infty([r^{-\vartheta}T+1/2,r^{-\vartheta}T+1]\times B_{1/2})} \leq C.
$$
By scaling, this implies that
$$
\|\partial_tu\|_{\LL^\infty([T+r^{\vartheta}/2,r^{-\vartheta}T+r^{-\vartheta}]\times B_{r}(x_0))} \leq Cr^{\frac{\gamma}{m}-\vartheta},
$$
and by the arbitrariness of $T$ and $x_0$ we obtain (recall that $r=\textrm{dist}(x_0,\partial\Omega)/2$)
\begin{equation}
\label{eq:ut}
|\partial_tu(t,x)|\leq C\textrm{dist}(x,\partial\Omega)^{\frac{\gamma}{m}-\vartheta}
\qquad \forall\, t \in [T_2,T_1],\,x \in \Omega.
 \end{equation}
 Note that $\frac{\gamma}{m}-\vartheta=\gamma-2s<0$ by our assumption.

 Now, given $t_0,t_1 \in [T_2,T_1]$ and $x \in \Omega$, we argue as follows:
 if $|t_0-t_1|\leq \textrm{dist}(x,\partial\Omega)^\vartheta$ then we
 use \eqref{eq:ut} to get (recall that $\frac{\gamma}{m}-\vartheta<0$)
 $$
 |u(t_1,x)-u(t_0,x)|\leq C\textrm{dist}(x,\partial\Omega)^{\frac{\gamma}{m}-\vartheta}|t_0-t_1|\leq C|t_0-t_1|^\frac{\gamma}{m\vartheta}.
 $$
 On the other hand, if $|t_0-t_1|\geq \textrm{dist}(x,\partial\Omega)^\vartheta$,
 then we use \eqref{eq:holder x} and the fact that $u$ vanishes on $\partial \Omega$
 to obtain
 $$
 |u(t_1,x)-u(t_0,x)| \leq |u(t_1,x)|+|u(t_0,x)| \leq C\textrm{dist}(x,\partial\Omega)^{\frac{\gamma}{m}}\leq C|t_0-t_1|^\frac{\gamma}{m\vartheta}.
 $$
 This proves that $u$ is $\frac{\gamma}{m\vartheta}$-H\"older continuous in time,
 and completes the proof of Theorem \ref{thm.regularity.2}. \qed
%
\section{Numerical evidence}\label{sec.numer}

After discovering the unexpected boundary behavior, we looked for  numerical confirmation.
This has been given to us by the authors of \cite{numerics}, who exploited the analytical tools developed in this paper to support our results by means of accurate numerical simulations. We include here some of these simulations, by courtesy of the authors. In all the figures we shall consider the Spectral Fractional Laplacian, so that $\gamma=1$ (see Section \ref{ssec.examples} for more details).

\begin{figure}[H]\label{fig1}
  \hspace{-1cm}
  \includegraphics[width=260pt]{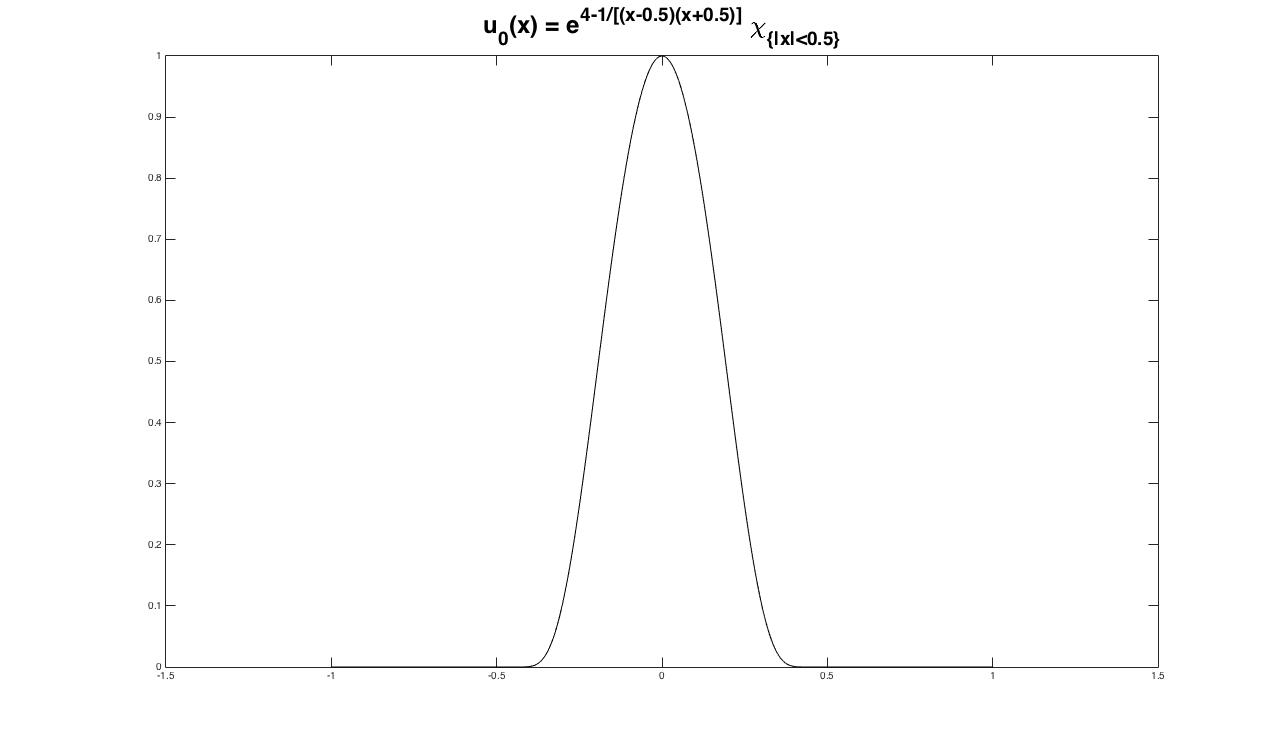}\hspace{-8mm}\includegraphics[width=260pt]{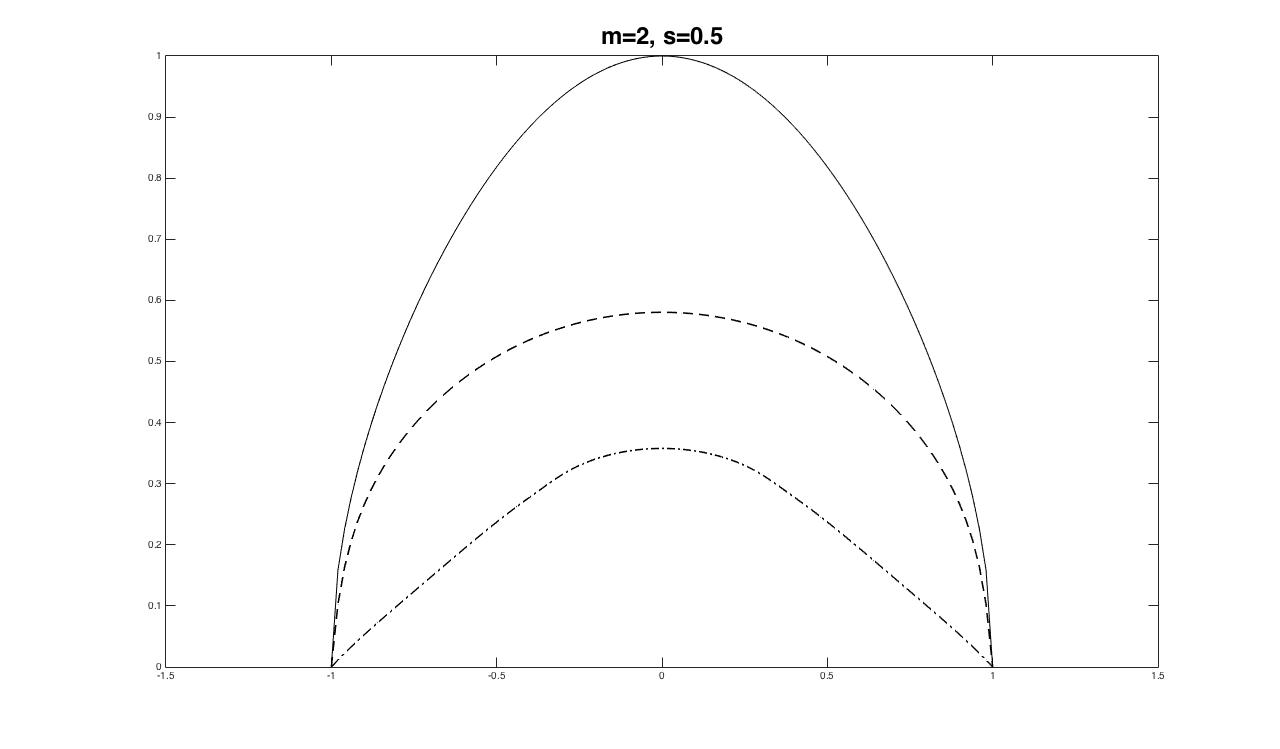}
  \vspace{-8mm}\caption{ \footnotesize
On the left, the initial condition $u_0$. On the right, the solid line represents $\Phi_1^{1/m}$, and the dotted lines represent $t^{\frac{1}{m-1}}u(t)$ at $t=1$
and $t=5$. The parameters are $m=2$ and $s=1/2$, hence $\sigma=1$.
While $u(t)$ appears to behave as $\Phi_1\asymp \dist(\cdot,\partial\Omega)$ for very short times, already at $t=5$ it exhibits the matching boundary behavior predicted by Theorem \ref{thm.Lower.PME.large.t}.
}
\end{figure}

We take $\Omega=(-1,1)$, and we consider as initial datum the compactly supported function $u_0(x)=e^{4-\frac{1}{(x-1/2)(x+1/2)}}\chi_{|x|< 1/2}$ appearing in the left of Figure \ref{fig1}.
In all the other figures, the solid line represents either \ $\Phi_1^{1/m}$ or $\Phi_1^{1-2s}$,
while the dotted lines represent \ $t^{\frac{1}{m-1}}u(t)$ for different values of $t$, where $u(t)$ is the solution starting from $u_0$.
These choices are motivated by Theorems \ref{thm.Lower.PME.large.t} and \ref{prop.counterex2}.
Since the map $t\mapsto t^{\frac{1}{m-1}}\,u(t,x)$ is nondecreasing for all $x \in \Omega$
(cf. (2.3) in \cite{BV-PPR2-1}),
the lower dotted line corresponds to an earlier time with respect to the higher one.
\vspace{-3mm}
\vspace{-4mm}

\begin{figure}[H]\label{fig2}
  \hspace{-1.2cm}
  \includegraphics[width=260pt]{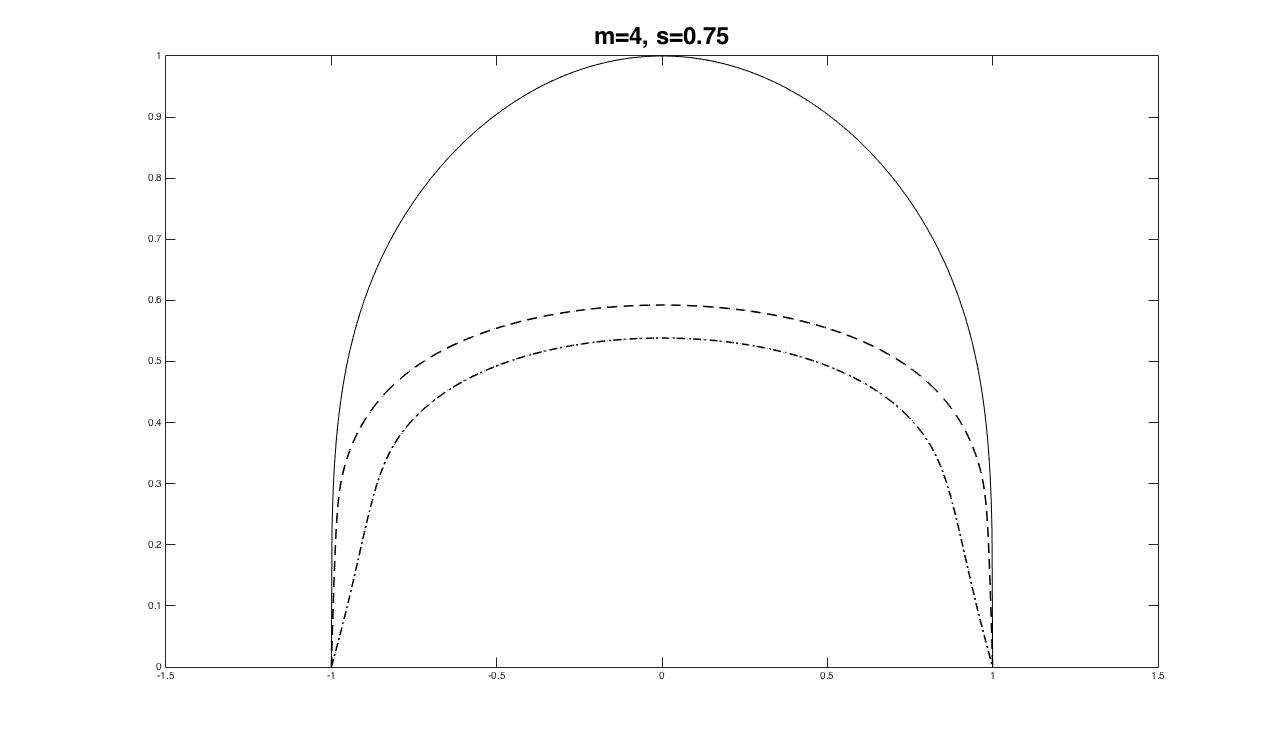}\hspace{-7mm} \includegraphics[width=260pt]{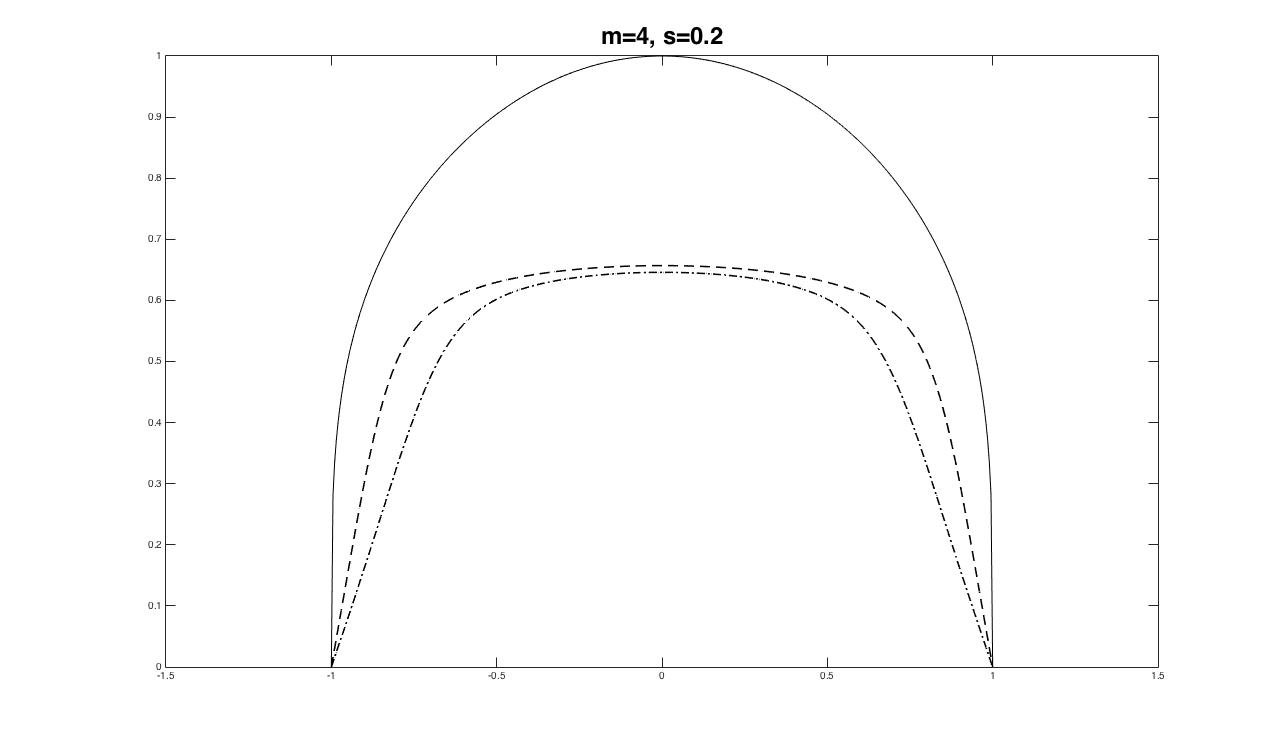}\vspace{-4mm}
   \caption{ \footnotesize
   In both pictures,  the solid line represents $\Phi_1^{1/m}$.
On the left, the dotted lines represent $t^{\frac{1}{m-1}}u(t)$ at $t=30$
and $t=150$, with parameters $m=4$ and $s=3/4$ (hence $\sigma=1$).
In this case $u(t)$ appears to behave as $\Phi_1\asymp \dist(\cdot,\partial\Omega)$ for quite some time, and only around $t=150$ it exhibits the matching boundary behavior predicted by Theorem \ref{thm.Lower.PME.large.t}.
On the right, the dotted lines represent $t^{\frac{1}{m-1}}u(t)$ at $t=150$ and $t=600$ with parameters $m=4$ and $s=1/5$ (hence $\sigma=8/15<1$).
In this case $u(t)$ seems to exhibit a linear boundary behavior even after long time (this linear boundary behavior is a universal lower bound for all times by Theorem
\ref{thm.Lower.PME}).
The second picture may lead one to conjecture that, in the case $\sigma<1$ and $u_0\lesssim \Phi_1$, the behavior $u(t)\asymp \Phi_1$ holds for all times.
However, as shown in Figure 3, there are cases when $u(t)\gg \Phi_1^{1-2s}$ for large times.}
\end{figure}\vspace{-4mm}

\begin{figure}[H]\label{fig3}
  \hspace{-1.2cm}
  \includegraphics[width=260pt]{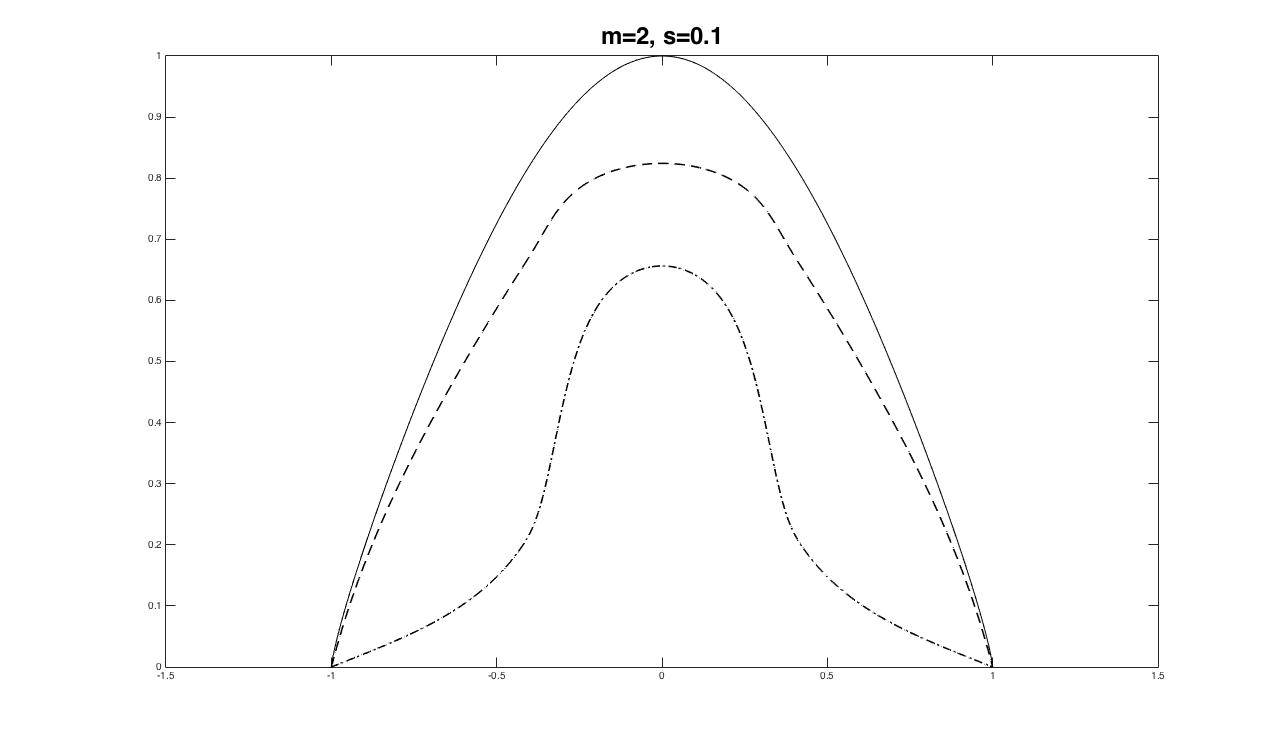}\hspace{-7mm}\includegraphics[width=260pt]{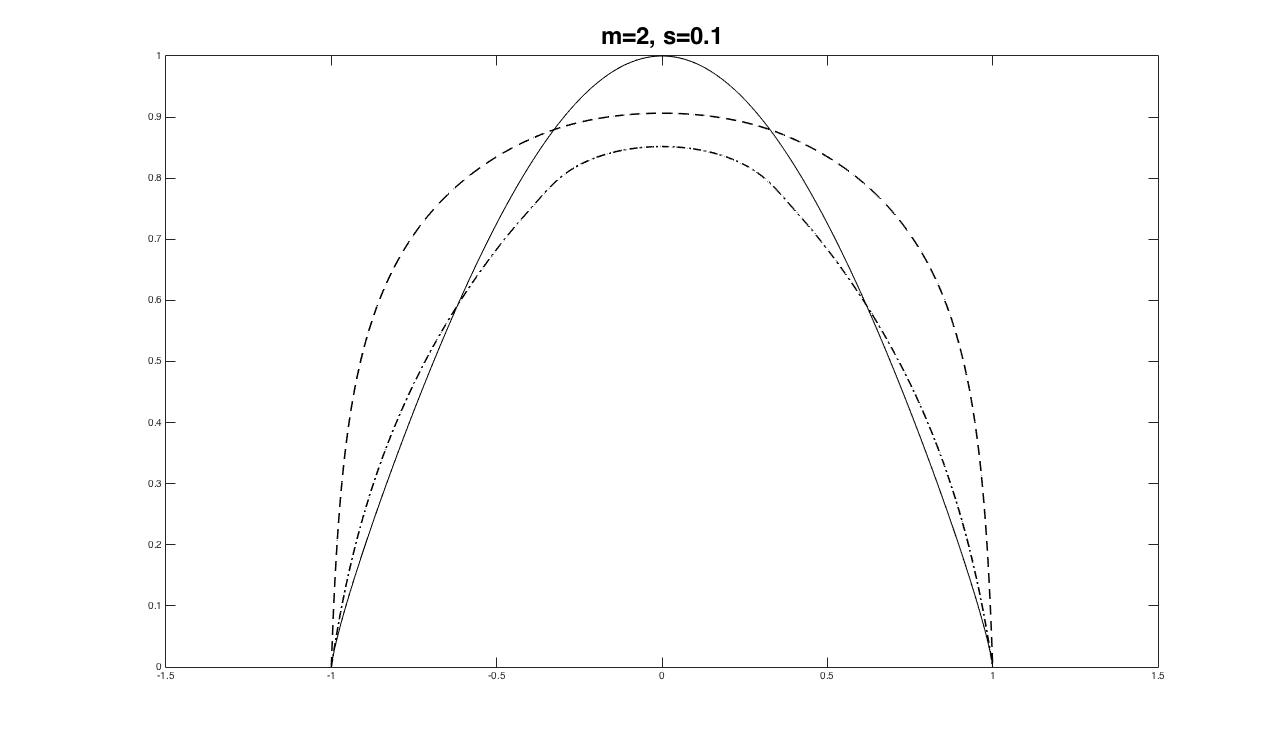}
  \vspace{-8mm}
  \caption{ \footnotesize
  In both pictures we use the parameters $m=2$ and $s=1/10$ (hence $\sigma=2/5<1$), and  the solid line represents $\Phi_1^{1-2s}$.
On the left, the  dotted lines represent $t^{\frac{1}{m-1}}u(t)$ at $t=4$
and $t=25$, on the right we see $t=40$ and $t=150$.
Note that $u(t)\asymp \Phi_1$ for short times.
Then, after some time, $u(t)$ starts looking more like $\Phi_1^{1-2s}$, and
for large times ($t=150$) it becomes much larger than  $\Phi_1^{1-2s}$.
 }
\end{figure}

Comparing Figures 2 and 3, it seems that when $\sigma<1$ there is no hope to find a universal behavior of solutions for large times. In particular, the bound provided by \eqref{intro.1} seems to be optimal.

\section{Complements, extensions and further examples}\label{sec.comm}

\noindent {\bf Elliptic versus parabolic.} The exceptional boundary behaviors we have found for some operators and data came as a surprise to us, since the solution to the corresponding ``elliptic setting'' $\A S^m= S$ satisfies $S\asymp \Phi_1^{\sigma/m}$  (with a logarithmic correction when
$2sm\ne \gamma(m-1)$),  hence separate-variable solutions always satisfy  \eqref{intro.1b} (see \eqref{friendly.giant} and Theorem \ref{Thm.Elliptic.Harnack.m}).

\medskip

\noindent\textbf{About the kernel of operators of the spectral type. }\label{ss2.2}
In this paragraph we study the properties of the kernel of $\A$. While in some situations $\A$ may not have a kernel (for instance, in the local case), in other situations that may not be so obvious from its definition. In the next lemma it is shown in particular that the SFL, defined by \eqref{sLapl.Omega.Spectral}, admits representation of the form \eqref{SFL.Kernel}. We state hereby the precise result, mentioned in \cite{AbTh} and proven in \cite{SV2003} for the SFL.
\begin{lem}[Spectral Kernels]\label{Lem.Spec.Ker}
Let $s \in (0,1)$, and let $\A$ be the $s^{th}$-spectral power of a linear elliptic second order operator $\mathcal A$, and let $\Phi_1 \asymp \dist(\cdot,\partial\Omega)^\gamma$ be the first positive eigenfunction of $\mathcal A$. Let $H(t,x,y)$ be the Heat Kernel of $\mathcal A$ , and assume that it satisfies the following bounds: there exist constants $c_0,c_1,c_2>0$ such that
for all $0<t\le1$
\begin{equation}\label{Bounds.HK.s=1.t<1}
 c_0\left[\frac{\Phi_1(x)}{t^{\gamma/2}}\wedge 1\right] \left[\frac{\Phi_1(y)}{t^{\gamma/2}}\wedge 1\right]\frac{\ee^{-c_1\frac{|x-y|^2}{t}}}{t^{N/2}}\le H(t,x,y)\le c_0^{-1}\left[\frac{\Phi_1(x)}{t^{\gamma/2}}\wedge 1\right] \left[\frac{\Phi_1(y)}{t^{\gamma/2}}\wedge 1\right]\frac{\ee^{-\frac{|x-y|^2}{c_1\,t}}}{t^{N/2}}
\end{equation}
and
\begin{equation}\label{Bounds.HK.s=1.t>1}
0\leq H(t,x,y) \leq c_2 \Phi_1(x)\Phi_1(y)\qquad\mbox{for all $t\ge 1$. }
\end{equation}
Then the operator $\A$ can be expressed in the form
\begin{equation}\label{A.op.kern.1}
\A f(x)=P.V.\int_{\RR^N} \big(f(x)-f(y)\big)\,K(x,y)\dy + B(x)u(x)
\end{equation}
with a kernel $K(x,y)$ supported in $\overline{\Omega}\times\overline{\Omega}$ satisfying
\begin{multline}\label{Lem.Spec.Ker.bounds}
K(x,y)\asymp\frac{1}{|x-y|^{N+2s}}
\left(\frac{\Phi_1(x)}{|x-y|^\gamma }\wedge 1\right)
\left(\frac{\Phi_1(y)}{|x-y|^\gamma }\wedge 1\right)\quad\mbox{and}\quad B(x)\asymp \Phi_1(x)^{-\frac{2s}{\gamma}}
\end{multline}
\end{lem}
The proof of this Lemma follows the ideas of \cite{SV2003}; indeed assumptions of Lemma \ref{Lem.Spec.Ker} allow to adapt the proof of \cite{SV2003} to our case with minor changes.

\smallskip

\noindent {\bf Method and generality. }Our work is part of a current effort aimed at extending the theory of evolution equations of parabolic type to a wide class of nonlocal operators, in particular operators with general kernels  that have been studied by various authors (see for instance \cite{EJT, DPQR, Ser14}).
Our approach is different from many others: indeed, even if the equation is nonlinear,  we concentrate on the properties of the inverse operator $\AI$ (more precisely, on its kernel given by the Green function $\K$), rather than on the operator $\A $ itself. Once this setting is well-established and  good linear estimates for the Green function are available, the calculations and estimates are very general. Hence, the method is applicable to a very large class of equations, both for Elliptic and Parabolic problems, as well as to more general nonlinearities than $F(u)=u^m$ (see also related comments in the works \cite{BV-PPR1, BSV2013, BV-PPR2-1}).

\smallskip

\noindent {\bf Finite and infinite propagation.} In all cases consider in the paper for $s<1$ we prove that the solution becomes strictly positive inside the domain at all positive times. This is called {\sl infinite speed of propagation}, a property that does not hold in the limit $s=1$ for any $m>1$ \cite{VazBook} (in that case, finite speed of propagation holds and a free boundary appears). Previous results on this infinite speed of propagation can be found in  \cite{BFR, DPQRV2}. We recall that infinite speed of propagation is typical of the evolution with nonlocal operators representing long-range interactions, but it is not true for the standard porous medium equation, hence a trade-off takes place when both effects are combined; all our models fall on the side of infinite propagation, but we recall that finite propagation holds for a related nonlocal model called ``nonlinear porous medium flow with fractional potential pressure'', cf. \cite{CV1}.

\smallskip

\noindent {\bf The local case.}  Since $2sm>\gamma(m-1)$ when $s=1$ (independently of $m>1$), our results give a sharp behavior in the local case after a ``waiting time''.
Although this is well-known for the classical porous medium equation, our results apply also to the case uniformly elliptic operator in divergence form with $C^1$ coefficients, and yield new results in this setting. Actually one can check that, even when the coefficients are merely measurable, many of our results are still true and they provided universal upper and lower
estimates.
At least to our knowledge, such general results are
completely new.

\subsection{Further examples of operators}\label{sec.examples}

Here we briefly exhibit a number of examples to which our theory applies, besides the RFL, CFL and SFL already discussed in Section \ref{sec.hyp.L}. These include a wide class of local and nonlocal operators. We just sketch the essential points, referring to \cite{BV-PPR2-1} for a more detailed exposition.

\medskip

\noindent\textbf{Censored Fractional Laplacian (CFL) and operators with more general kernels. }As already mentioned in Section \ref{ssec.examples}, assumptions (A1), (A2), and (K2) are satisfied  with $\gamma=s-1/2$. Moreover, it follows by \cite{bogdan-censor, Song-coeff} that we can also consider operators of the form:
\[
\A f(x)=\mathrm{P.V.}\int_{\Omega}\left(f(x)-f(y)\right)\frac{a(x,y)}{|x-y|^{N+2s}}\dy\,,\qquad\mbox{with }\frac{1}{2}<s<1\,,
\]
where $a(x,y)$ is a symmetric function of class $C^1$ bounded between two positive constants. The Green function $\K(x,y)$ of $\A$ satisfies the stronger assumption (K4), cf. Corollary 1.2 of~\cite{Song-coeff}.

\medskip

\noindent\textbf{Fractional operators with more general kernels. }Consider integral operators of the form\vspace{-1mm}
\[
\A f(x)=\mathrm{P.V.}\int_{\RR^N}\left(f(x)-f(y)\right)\frac{a(x,y)}{|x-y|^{N+2s}}\dy\,,
\]
where $a$ is a measurable symmetric function, bounded between two positive constants, and satisfying
\[
\big|a(x,y)-a(x,x)\big|\,\chi_{|x-y|<1}\le c |x-y|^\sigma\,,\qquad\mbox{with }0<s<\sigma\le 1\,,
\]
for some  $c>0$ (actually, one can allow even more general kernels, cf. \cite{BV-PPR2-1, Kim-Coeff}).
 Then, for all $s\in (0, 1]$, the Green function $\K(x,y)$ of $\A$ satisfies (K4) with $\gamma=s$\,, cf. Corollary 1.4 of \cite{Kim-Coeff}.

\medskip

\noindent\textbf{Spectral powers of uniformly elliptic operators. }Consider a linear operator $\mathcal A$ in divergence form,
\[
\mathcal A=-\sum_{i,j=1}^N\partial_i(a_{ij}\partial_j)\,,
\]
with uniformly elliptic $C^1$ coefficients. The uniform ellipticity allows one to build a self-adjoint operator on $\LL^2(\Omega)$ with discrete spectrum $(\lambda_k, \phi_k)$\,. Using the spectral theorem, we can construct the spectral power of such operator as follows
\[
\A f(x):=\mathcal A^s\,f(x):=\sum_{k=1}^\infty \lambda_k^s \hat{f}_k \phi_k(x),\qquad\mbox{where }\qquad \hat{f}_k=\int_\Omega f(x)\phi_k(x)\dx
\]
(we refer to the books \cite{Davies1,Davies2} for further details), and the Green function satisfies (K2) with $\gamma=1$\,, cf. \cite[Chapter 4.6]{Davies2}. Then, the first eigenfunction $\Phi_1$
is comparable to $\dist(\cdot, \partial\Omega)$.
Also, Lemma \ref{Lem.Spec.Ker} applies (see for instance \cite{Davies2}) and allow us to get sharp upper and lower estimates for the kernel $K$ of $\A$, as in \eqref{Lem.Spec.Ker.bounds}\,. 

\noindent\textbf{Other examples. } As explained in Section 3 of \cite{BV-PPR2-1}, our theory may also be applied to: (i) Sums of two fractional operators; (ii) Sum of the Laplacian and a nonlocal operator kernels; (iii) Schr\"odinger equations for non-symmetric diffusions; (iv) Gradient perturbation of restricted fractional Laplacians.
Finally, it is worth mentioning that our arguments readily extend to operators on manifolds
for which the required bounds hold.

\vskip .5cm

\noindent{\bf Acknowledgments. }M.B. and J.L.V. are partially funded by Project MTM2011-24696 and  MTM2014-52240-P(Spain). A.F. has been supported by NSF Grants DMS-1262411 and DMS-1361122, and by the ERC Grant ``Regularity and Stability in Partial Differential Equations (RSPDE)''. M.B. and J.L.V. would like to acknowledge the hospitality of the Mathematics Department of the University of Texas at Austin, where part of this work has been done. J.L.V. was also invited by BCAM, Bilbao. We thank an anonymous referee for pointing out that Lemma \ref{Lem.Spec.Ker} was proved in  \cite{SV2003} and mentioned in \cite{AbTh}.

\end{document}